\documentclass[leqno]{article}
\usepackage{latexsym,amsmath,amsfonts,mathrsfs,amssymb,amsthm}
\usepackage{yhmath}
\usepackage[usenames]{color}
\usepackage{graphicx}
\def\R{\mathbb R}
\def\N{\mathbb N}

\def\Z{\mathbb Z}
\def\V{\mathbb V}
\def\s{\sharp}
\def\p{\partial}
\def\a{\alpha}
\def\e{\epsilon}
\def\d{\delta}
\def\H{{\cal H} }
\def\be{\begin{equation}}
\def\ee{\end{equation}}
\def\bs{\backslash}
\def\qed{\hfill$\Box$\bigskip}

\def\nd{\noindent \textbf{Proof.} }
\numberwithin{equation}{section}
\newtheorem{lem}[equation]{Lemma}
\newtheorem{pro}[equation]{Proposition}
\newtheorem{defn}[equation]{Definition}
\newtheorem{ex}[equation]{Example}
\newtheorem{thm}[equation]{Theorem}

\newtheorem{rem}[equation]{Remark}
\begin{document}
\bigskip

\begin{center}\Large \textbf{Unique tangent behavior for 1-dimensional stationary varifolds}\end{center}

\bigskip

\centerline{\large Xiangyu Liang}

%

\vskip 1cm

\centerline {\large\textbf{Abstract.}}

We prove that, without any assumption on lower density bound or codimension, any 1-dimensional stationary varifold on any Riemannian manifold admits unique tangent behaviour everywhere.

\bigskip


\bigskip


\setcounter{section}{-1}

\section{Introduction}

In geometric measure theory, one main interest is the theory of minimal sets, currents and varifolds, which aims at understanding the existence and regularity of physical objects that admit certain minimizing property.  This is known as Plateau's problem in physics.

 Lots of notions of minimality have been introduced to modernise this problem, such as minimal surfaces, mass or size minimizing integral currents (De Rham, Federer-Fleming), minimal sets (Almgren, Reifenberg) and stationary varifolds (Almgren).
 
 Unlike the traditional minimal surface theory, these are all objects with singularities. In order to understand the local structure around singular points, a first step is to study the tangent behavior. Here, the notion of tangent (see Definition 1.18 for the case of varifolds for example) is in a much weaker sense than the tangent space in differential geometry, such that one can have more than one tangent objects at a point even for nice objects. For example in  \cite{Kor11,Kor11a}, the authors construct a set with rather good regularity: almost smooth, with monotone density ratio; but still, it admits more than one tangent cone at a point.

As a result, in most circumstances, the uniqueness of tangent behavior is always an important regularity property, and it has been widely investigated in many of the above theories, as well as in other problems in analysis, e.g. currents (\cite{Bel14, Whi83, PuRi10, Kis88}), harmonic maps (\cite{RiTi04,Whi92}), minimal sets (\cite{DEpi, Ta}), stationary varifolds (\cite{AlAl,HP09}), minimal surfaces (\cite{Sim83a, AlAl81}), etc.. 

In this article we discuss the unique tangent behavior for stationary 1-varifolds. Stationary varifolds are weak solutions for Plateaus problem in the setting of measures, defined as critical points of measure while deforming along any vector fields. See Definition 1.23. The notion of stationary varifold is perhaps one of the weakest: minimal surfaces, minimal sets (in the sense of Almgren or Reifenberg), size minimizing integral currents, etc., are all stationary varifolds.

The question about uniqueness of tangent behaviour for stationary varifolds was raised by Allard (\cite{All72} 6.5), and then by Simon (\cite{Sim83} 42.2). A number of partial results are known: for non regular case, there is a counter example constructed in \cite{Kor15}: there exists a stationary 2-rectifiable varifold (but without lower density bounds, hence the support of the varifold is 3 dimensional) which admits more than one tangent at a point with positive density; on the other hand, there are a number of cases where we have affirmative answers: in \cite{AlAl} the authors prove that any 1-dimensional stationary varifolds with a uniform lower density bound (in particular they are rectifiable) in any Riemannian manifold  admits a unique tangent varifold at every point; and in \cite{HP09}, the author proves that in $\R^2$, any stationary 1-varifold admit unique tangent behaviour at every point.

While the picture for a general answer to arbitrary dimensional and codimensional cases is still obscure, here we give an affirmative answer to any stationary 1-varifold in any Riemannian manifold:

\begin{thm} Let $M$ be an $n$-dimensional Riemannian manifold, and let $V$ be a stationary 1-varifold on $M$. Then for any $x_0\in M$ where $V$ has non zero density, $V$ admits only one tangent varifold at $x_0$. In particular, if $V$ is a stationary 1-varifold in an open set $U\subset\R^n$, then at every point where $V$ has non zero density, $V$ admits only one tangent varifold.
\end{thm}

As for the proof, a traditional way in analysis of proving uniqueness of tangent behavior is to show that the decay of the density function is fast enough, and this often envolves the proof of an epiperimetric property. See for example (\cite{AlAl,Whi83,PuRi10, RiTi09,Sim83a,DEpi, Ta}). Alternative proofs have been given by calibrations (see \cite{Bel11,Bel14}), or a representation in a homological way (see \cite{HP09}).

Our idea of proof is somehow different from all the above. On the other hand, the result of \cite{HP09}, which is a particular case of Theorem 0.1, is in fact the base of our proof. So let us say a little more:

In \cite{HP09}, the author proves that every 1-dimensional stationary varifold is in some sense the second derivative of a convex function. Then the unique tangent behavior follows directly from this representation and properties of convex functions. However, the representation is no longer true even in dimension 3 (as remarked at the end of \cite{HP09}), and up to now we do not find an easy generalization of this representation.

Nevertheless, although we fail to apply the method in \cite{HP09} to our case, the result in \cite{HP09} itself is quite useful for us. We can see this  in the following sketch of proof of Theorem 0.1:


We first prove in Section 2 that weighted projections (different from mapping varifolds under projections, see definition 2.1) of any stationary 1-varifold are still stationary 1-varifolds, provided that the projections are locally finite for this varifold. We will then get rid of the hypothesis of projections being locally finite, using some cut and paste operation. In fact, one can notice that the restriction of a varifold to any bounded set is finite, and hence always admits locally finite projections. Therefore if the projection of the whole varifold is locally infinite, this must be caused by its behavior at infinity. On the contrary, the tangent behavior of a varifold at any point $x$ is a local behavior. Hence we use the cut and paste operation to get rid of the part at infinity: we keep the part of a stationary 1-varifold $V$ in the unit ball $B_x$ centered at the point $x$, throw away the outer part, and paste some other simpler varifold $V'$ (supported outside the ball $B_x$) to get a new stationary varifold $W=V\lfloor_{B_x}+V'$. The definition of $V'$ will depend on the first variation of $V\lfloor_{B_x}$, and it will admit locally finite projections. The varifold $W$ admits same tangent varifolds as $V$ at $x$. See Section 3 for details.

As a consequence, since we know that stationary 1-varifolds in all 2-dimensional subspaces admit unique tangent varifold at every point, the question turns to be that, if two tangent varifolds (which are conic) of dimension 1 admit the same weighted projection to every direction, are they the same ?

%

In our proof, we discuss two cases: in rectifiable case, we prove that same projections on $n-1$-subspaces in a set of positive measure in the Grassmannian will guarantee the same tangent varifold (Section 4); for the general case, we do not know whether the same projections on a subset of Grassmannian still works, so we need the same projection on almost all directions (see Section 5). 

As a remark, notice that the above result is not true in dimension 2 for general 1-varifolds. See Remark 4.61. For rectifiable case, our proof uses the fact that the ambient dimension is at least 3 (Still see Remark 4.62). But we do not know whether in dimension 2, weighted projections can determine a conic rectifiable 1-varifold.

The last section is devoted to the case for a general Riemannian manifold $M$. The proof is based on Nash embedding theorem and the cut and paste procedure as before.

\textbf{Acknowledgement:} The author would like to thank David Preiss for introducing this problem and the previous work \cite{HP09}, and thank him and Guy David for many helpful discussions.

\section{Definitions and preliminaries}

\subsection{Basic notations}

Let $a,b\in \R^n$.

$[a,b]$ denotes the segment with endpoints $a$ and $b$;

$\overrightarrow{ab}$ is the vector $b-a$;

$S^{n-1}$ denotes the unit sphere in $\R^n$;

For any $z\in \R^n\bs \{0\}$, $R_z$ denotes the half line issued from the origin and passing through $z$;

For any subspaces $P$ of $\R^n$, and any $x\in P$, let $B_P(x,r):=\{y\in P, |y-x|<r\}$ be the open ball in $P$, centered at $x$ with radius $r$. Denote by $\pi_P:\R^n\to P$ the orthogonal projection. And for any $x\in \R^n$, let $x_P=\pi_P(x)\in P$.

For any two different points $a,b\in\R^n$, $R_{a,b}$ denotes the half line issued from $a$ and passing through $b$. For any point $a$ different from the origin, $R_a$ denotes the half line issued from the origin and passing through the point $a$;

For $k\le n$, $\H^k$ denotes the $k$-dimensional Hausdorff measure.

\subsection{Varifolds}
 
Here all the notations concerning varifolds are mainly those in \cite{All72} and \cite{Sim83}. All sets in all the definitions below are measurable. Although in our paper we only discuss varifolds of dimension 1, definitions and theorems are given in general dimensions, since they are not more complicated in general dimensions.

\begin{defn}[Grassmann Algebra. cf. \cite{Fe} Chapter 1] Let $V$ be a finite dimensional Euclidean vector space, Let $n=dim V$,

$1^\circ$ For $k\le m$, let $\wedge_kV$ be the $k$-th exterior algebra of $V$. Let $G(V,k)$ be the space of $k$ dimensional linear subspaces of $V$. For any $S\in G(V,k)$, set $S^\perp=\{x:x\cdot y=0\mbox{ for all }y\in S\}\in G(V,n-k)$, and let $\pi_S$ denote the orthogonal projection from $V$ to $S$. 

In particular, when $V=\R^n$, let $G(n,k)$ denote the space of $k$ dimensional linear subspaces of $\R^n$; 

$2^\circ$ Let $f: \R^n\to \R^n$ be a linear map, and $S\in G(n,k)$. Then let $|f\circ\pi_S|$ denotes $|\wedge_k(f\circ \pi_S)|$. And let $f(S)$ denotes the image $\{f(x);x\in \R^n\}$. When $|f\circ\pi_S|\ne 0$, $f(S)\in G(n,k)$;

$3^\circ$ Let $P\in G(n,m)$. For any $k\le m$, set $G_k(P)=P\times G(P,k)$. In particular, denote by $G_k(n)$ the product $\R^n\times G(n,k)$. 
\end{defn}

\begin{defn}[Varifold]Let $k\le n$.

$1^\circ$ A $k$-dimensional varifold $V$ in $\R^n$ is a Radon measure on $G_k(n)$. Let $\V_k(\R^n)$ denote the space of $k$-dimensional varifolds in $\R^n$;

$2^\circ$ For $V\in \V_k(\R^n)$, the total mass $|V|$ of $V$ is the Radon measure on $\R^n$ defined by 
\be |V|(A)=V(A\times G(n,k)),\mbox{ for all }A\subset \R^n.\ee

We say that a $k$-varifold $V$ is finite if $|V|$ is a finite measure.
\end{defn}

\begin{defn}[Rectifiable varifolds]Let $E\subset \R^n$ be a $k$-dimensional rectifiable set with locally finite Hausdorff measure. For $\H^k$-a.e. $x\in E$, denote by $T_xE$ the $k$-dimensional approximate tangent plane of $E$ at $x$. Let $V_E$ denote the varifold defined as
\be V_E(A)=\H^k\{x:(x,T_xE)\in A\},\forall A\in G_k(n),\ee
and called the $k$-varifold induced by $E$.
Given a positive $\H^k\lfloor_E$-integrable function $\theta$, define $V_{E,\theta}$ as 
\be V_{E,\theta}(A)=\int_E \theta(x)1_A(x,T_xE)d\H^k, \forall A\in G_k(n).\ee
Clearly we have
\be |V_E|=\H^k\lfloor_E\ and\ |V_{E,\theta}|=\theta\H^k\lfloor_E,\ee
where $\theta\H^k\lfloor_E$ denotes the measure
\be\theta\H^k\lfloor_E(A)=\int_{A\cap E}\theta(x)d\H^k(x),\forall A\subset \R.\ee
We say a varifold $V\in V_k(\R^n)$ is $k$-dimensional rectifiable, if there exists a $k$-dimensional rectifiable set $E$ with locally finite Hausdorff measure, and a positive $\H^k\lfloor_E$-integrable function $\theta$ such that $V=V_{E,\theta}$. The function $\theta$ is called the multiplicity function, or the density function. When $\theta$ is $\H^k\lfloor_E$-a.e. an integer, then $V$ is called an integral varifold.
\end{defn}


\begin{defn}[Density for varifolds]$1^\circ$ Let $\mu$ be a Borel regular measure in $\R^n$, and $k$ be an integer. Whenever $x\in \R^n$, let 
\be\theta_{k*}(\mu,x)=\liminf_{r\to 0}\frac{\mu(B(x,r))}{\alpha(k)r^k},\ \theta^*_k(\mu,x)=\limsup_{r\to 0}\frac{\mu(B(x,r))}{\alpha(k)r^k},\ee
where $\a(k)$ denotes the $k$-Lebesgue measure of the unit ball of $\R^k$.
When the $\theta_{k*}(\mu,x)$ and $\theta^*_k(\mu,x)$ coincide, let
\be \theta_k(\mu,x)=\lim_{r\to 0}\frac{\mu(B(x,r))}{\alpha(k)r^k}.\ee
$\theta_{k*}$, $\theta^*_k$ and $\theta_k$ are called respectively the $k$-dimensional lower density, upper density and density of $\mu$ at $x$. For any $k$-rectifiable set $E$ with locally finite Hausdorff measure, the $k$-density of $\H^k\lfloor_E$ exists for $\H^k$ almost all points in the set.

$2^\circ$ Let $V\in \V_k(\R^n)$. Let $x\in \R^n$. Then the $k$-lower density, upper density and density $\theta_{k*}(V,x)$, $\theta^*_k(V,x)$ and $\theta_k(V,x)$ of $V$ at $x$ are just those of the measure $|V|$, that is,
\be \theta_{k*}(V,x):=\theta_{k*}(|V|,x), \theta^*_k(V,x):=\theta^*_k(|V|,x)\mbox{ and } \theta_k(V,x):=\theta_k(|V|,x).\ee

$3^\circ$ Let $V\in \V_k(\R^n)$. We denote by $[V]$ the set of points $x\in \R^n$ such that $\theta_k(|V|,x)$ exists and $0<\theta_k(|V|,x)<\infty$.
\end{defn}


\subsection{Mapping varifolds, tangents, and the first variation}

The next definitions will all rely on the image of a varifold by a smooth map, so let us first give the following definition:

\begin{defn}[Image by a smooth map] Let $f: \R^n\to \R^m$ be smooth. Let $V$ be a $k$-dimensional varifold in $\R^n$. Then the Borel regular measure $f_\s(V)$ on $G_k(\R^m)$, called the image of $V$ by $f$, is defined as follows: 
\be f_\s(V)(B)=\int_{\{(x,S):(f(x), Df(x)(S))\in B\}}|Df\circ\pi_S|dV(x,S), \forall \mbox{ Borel subset }B\mbox{ of }G_k(\R^m),\ee
or equivalently,
\be \int_{G_k(\R^m)} g(x,S) df_\s(V)(x,S)=\int_{G_k(\R^n)}g(f(x),Df(x)(S))|Df\circ\pi_S|dV(x,S)\ee
for all Borel function $g$ on $G_k(\R^m)$.

In particular, when $f_\s(V)$ is locally finite, we call it the image varifold of $V$ under $f$.
\end{defn}

\begin{rem}It is easy to see from Definition 1.13 that if $V=V_E$ is the $k$-varifold induced by a $k$-rectifiable set $E$, and $f_\s(V_E)$ is a varifold, then $f_s(V_E)=V_{f(E)}$, that is, the $k$-varifold induced by $f(E)$.
\end{rem}

Image varifold will be used in particular to define tangent varifolds, and stationary varifolds. Let us begin with the tangents. 

When talking about tangents, one always refer to asymptotic behavior at small scales. Indeed, we will use the following family of dilatation maps:
for $x\in \R^n$ and $\lambda>0$, set
\be\eta_{x,\lambda}(y)=\frac{y-x}{\lambda},y\in \R^n.\ee

\begin{defn}[Tangent Varifolds] Let $V, C\in \V_k(\R^n)$, $x\in \R^n$. We say that $C$ is a tangent varifold of $V$ at $x$, if there exists $\lambda_i>0$, and $\lambda_i\to 0$, such that 
\be {\eta_{x,\lambda_i}}_\s V\to C\ee
weakly. Or equivalently, for any $f\in C_c(G_k(n))$, 
\be \int f(y,S)dC(y,S)=\lim_{i\to\infty}(\lambda_i)^{-k}\int f(\eta_{x,\lambda_i}(y),S)dV(y,S).\ee
 Denote by Var Tan$(V,x)$ the set of all tangent varifold of $V$ at $x$.\end{defn}
 
 \begin{rem}1) It is easy to see that whenever $x\in [V]$, Var Tan$(V,x)$ is non empty and compact. In this case, for any $C\in$Var Tan$(V,x)$, 
 \be|C|(B(0,r))=r^k\a(k)\theta_k(V,x),\ee
  and in particular
 $\theta_k(C,0)=\theta_k(V,x)$.
 
 2) The case where $x\not\in [V]$ is trivial.
 \end{rem}

Next let us introduce the first variation, and stationary varifolds. 

\begin{defn}[First variation, stationary varifold] Let $V\in \V_k(\R^n)$. 

$1^\circ$ The first variation $\d V$ of $V$ is a linear map from the space of smooth compactly supported vector fields on $\R^n$ to $\R$:
\be \d V(g)=\int div_sg(x) dV(x,S),\forall g\in C_C^\infty(\R^n,\R^n),\ee
where $div_Xg(x)$ is the divergence of $g$ restricted to $S$, that is, if $e_1,\cdots e_k$ is an orthonormal basis of $S$, then $div_Sg(x)=\sum_{i=1}^k \frac{\partial (g\cdot e_i)}{\partial e_i}(x)$.

$2^\circ$ Denote by $|\d V|$ the total variation of $\d V$, which is the Borel regular measure on $\R^n$ determined by
\be |\d V|(B)=\sup \{\d V(g): g\in C_C^\infty(\R^n), spt\ g\subset B\mbox{ and }|g|\le 1\}\ee
for all Borel set $B$.

$3^\circ$ We say that $V$ is a stationary varifold in $\R^n$ if $\d V=0$. If $U\subset \R^n$ is open, and $\d V(g)=0$ for all $g\in C_C^\infty(U,\R^n)$ then we say that $V$ is stationary in $U$.
\end{defn}

The first variation of a varifold measures the rate of change of the mass of the varifold while deformed along each vector field $g$. If the varifold is defined by a $C^2$ $k$-submanifold, then the first variation coincides with the mean curvature. In this case, a varifold induced by a $C^2$ $k$-submanifold is stationary if the submanifold is a minimal surface. See  \cite{All72} or \cite{Sim83} for more details.

Stationary varifolds admit many good properties. The one that will be helpful for us is the following:

\begin{pro}[cf. \cite{Sim83} 40.5]Let $V\in \V_k(\R^n)$ be stationary in an open subset $U$ of $\R^n$. Then for all $x\in U$, the density $\theta_k(V,x)$ exists.
\end{pro}

For non-stationary varifolds, when $|\d V|$ is locally finite, we also have a good representation:

\begin{pro}[cf. \cite{All72} 4.3] Let $V\in \V_k(\R^n)$, and $|\d V|$ is a Radon measure. Then there exists a $|\d V|$ measurable function $\omega(V;\cdot)$ with values in $S^{n-1}$ such that
\be \d V(g)=\int <g(x)\cdot \omega(V,x)> d|\d V|(x), \forall g\in C_C^\infty(\R^n).\ee
Evidently, $\omega(V,\cdot)$ is $|\d V|$ almost unique.
\end{pro}

\section{Weighted projection of stationary 1-varifolds}

From now on, we will concentrate ourselves on 1-varifolds.

\begin{defn}[Weighted projection varifold]Let $V$ be a 1-dimensional varifold in $\R^n$, and let $P$ be any plane of codimension at least 1 in $\R^n$. Let $\pi_P$ denote the orthogonal projection from $\R^n$ to $P$. We define the weighted projection varifold ${\pi_P}_{**}(V)$ of $V$ to $P$ to be the Borel regular measure on $G_1(P)$ defined as follows:
\be {\pi_P}_{**}(V)(B)=\int_{\{(x,S):(\pi_P(x), \pi_P(S))\in B\}}|\pi_P\circ \pi_S|^2dV(x,S)\ee
\end{defn}

\begin{rem} $1^\circ$ Note that this definition does not coincide with the usual mapping varifold ${\pi_P}_\s(V)$ of $V$ under the mapping $\pi_P$ (which is more often used, See Definition 1.13). For example let $L\in \R^2$ be the line generated by the vector $e_1+e_2$, and $V_L$ denote the rectifiable 1-varifold associate to $L$. Let $P$ be the line generated by $e_1$. Then ${\pi_P}_\s(V)$ is just the rectifiable 1-varifold $V_P$ associated to $P$, but ${\pi_P}_{**}(V)$ is $\frac{\sqrt 2}{2}V_P$.

$2^\circ$ We can easily check that $\forall B\subset P$, $|{\pi_P}_{**}(V)|(B)\le |{\pi_P}_\s(V)|(B)\le|V|(B\times P^\perp)$. As a result, if $|V|$ is finite, then ${\pi_P}_\s(V)$ and ${\pi_P}_{**}(V)$ are both finite, and hence are both 1-varifolds.
\end{rem}

In the rest of the paper, we will call ${\pi_P}_\s(V)$ the projection image of $V$, and ${\pi_P}_{**}(V)$ the weighted projection of $V$. For the weighted projection, we have the following similar property as projection images.

\begin{lem}Let $Q\subset P$ be two linear subspaces of $\R^n$. Then for any $V\in \V_1(\R^n)$,
\be {\pi_Q}_{**}(V)={\pi_Q}_{**}\circ {\pi_P}_{**}(V).\ee
\end{lem}

\nd This comes directly from the definition. \qed

\begin{pro} Let $V$ be a 1-dimensional stationary varifold in $\R^n$. Then for any $d<n$ and any $d$-plane $P$ in $\R^n$, if ${\pi_P}_\s(V)$ is locally finite, then ${\pi_P}_{**}(V)$ is a stationary varifold in $P$.
\end{pro}

\nd Without loss of generality, we can suppose that $P$ is the plane $\{(x_1,\cdots, x_n):x_{d+1}=x_{d+2}=\cdots =x_n=0\}$. Let $g$ be any smooth vector field with compact support on $P$, we would like to show that
\be \d {\pi_P}_{**}(V)(g)= \int_{G_1(P)}div_S g(x)d{\pi_P}_{**}(V)(x,S)=0.\ee

Now for any point $z\in \R^n$, we write $z=(x,y)$ where $x\in \R^d$ is the first $d$ coordinates of $z$, and $y\in \R^{n-d}$ the following coordinates of $z$. Thus $P=\{(x,y)\in \R^n:y=0\}$ and $P^\perp=\{(x,y)\in \R^n:x=0\}$.

Set $f:\R^n\to \R^n, f(x,y)=g(x)\in P$ for all $(x,y)\in\R^n$. Let $\chi:\R\to \R$ be a smooth lump function that satisfies:
\be \chi|_{[-\frac14, \frac14]}=1, \chi|_{(-\infty, -\frac34]\cup [\frac34,\infty)}=0, \chi(r)=\chi(-r), \mbox{ and }\chi\mbox{ is non increasing on }[0, 1].\ee

For each $m$, let $\chi_m:[0,\infty)\to \R$ be the smooth lump function:
\be \chi_m(r)=\left\{\begin{array}{cc}
1, &r\le m-1;\\
0, &r\ge m;\\
\chi(r-(m-1)), &m-1\le r\le m;\\
\chi(r+(m-1)), &-m\le r\le -m+1.\end{array}\right.\ee
Set $f_m(x,y)=f(x,y)\chi_m(|y|)$. Then $f_m$ are smooth vector fields with compact supports in $\R^n$, with $supp\ f\subset C_m$, where $C_m=B\times B_{\R^{n-d}}(0,m)$, and $B\subset P$ is a compact ball such that $supp\ g\subset B$. Moreover, $f_m=f$ on $C_{m-1}$.

Since $V$ is stationary in $\R^n$, for each $m$:
\be 0=\d V(f_m)=\int div_Sf_m(z)dV(z,S).\ee

Now for any $S\in G_1(\R^n)$, let $s$ be a unit vector that generates $S$. Since $\R^n=P\oplus P^\perp$, there is a unique decomposition $s=s_p+s^\perp$, with $s^\perp\in P^\perp$, and $s_p\in P$. Note that $|\pi_P\circ \pi_S|=|s_p|$. For any smooth vector field $h:\R^n\to \R^n$, 
\be div_S h(x)=\partial_s (h\cdot s)=\partial_{s^\perp+s_p}(h\cdot (s^\perp+s_p))=\p_{s^\perp}(h\cdot s^\perp)+\p_{s^\perp}(h\cdot s_p)+\p_{s_p}(h\cdot s^\perp)+\p_{s_p}(h\cdot s_p).\ee

For our maps $f_m$, their images are all in $P$, hence $f_m\cdot s^\perp=0$. Therefore,
\be div_S f_m=\p_{s^\perp}(f_m\cdot s_p)+\p_{s_p}(f_m\cdot s_p).\ee

Also, $\p_{s^\perp}(f_m\cdot s_p)=\p_{s^\perp}f=0$ for $(x,y)\in \R^n$ with $|y|\le m-1$, and $\p_{s^\perp} (f_m\cdot s_p)=\p_{s^\perp} (f\chi_m(|y|)\cdot s_p)=(f\cdot s_p)\p_{s^\perp}\chi_m(|y|)$ for $(x,y)\in \R^n$ with $|y|\in [m-1,m]$. 
Therefore, for $z=(x,y)\in C_{m-1} \mbox{ and }S\in G(n,1)$,
\be \begin{split}div_S f_m(x,y)&=\p_{s_p}(f_m\cdot s_p)=\p_{s_p}(f\cdot s_p)=|s_p|^2 div_{\pi_P(S)}f_m(x,y)\\
&=|s_p|^2 div_{\pi_P(S)}f(x,y)=|s_p|^2 div_{\pi_P(S)}g(x),\end{split}\ee
and for $z=(x,y)\in C_m\bs C_{m-1} \mbox{ and }S\in G(n,1)$,
\be \begin{split}
div_S f_m(x,y)&=(f\cdot s_p)\p_{s^\perp}\chi_m(|y|)+\p_{s_p}(f\cdot s_p)(x,y)\chi_m(|y|)\\
&=(g\cdot s_p)(x)\p_{s^\perp}\chi_m(|y|)+|s_p|^2 div_{\pi_P(S)}g(x)\chi_m(|y|).
\end{split}\ee

Of course, $div_Sf_m(z)=0$ outside $C_m$.
Hence by (2.10), 
\be \begin{split}0=&\int div_Sf_m(z)dV(z,S)=\int_{G_1(C_{m-1})}div_Sf_m(z)dV(z,S)+\int_{G_1(C_m\bs C_{m-1})}div_Sf_m(z)dV(z,S)\\
=&\int_{G_1(C_m\bs C_{m-1})}\{(g\cdot s_p)(x)\p_{s^\perp}\chi_m(|y|)+|s_p|^2 div_{\pi_P(S)}g(x)\chi_m(|y|)\}dV(z,S)\\
&+\int_{G_1(C_{m-1})}|s_p|^2 div_{\pi_P(S)}g(x)dV(z,S).\end{split}\ee

For the first term, note that $\chi_m$, $\chi_m'$, $g$ and $div g$ are bounded, and $\p_{s^\perp}\chi_m(|y|)\le |s^\perp|||\chi_m'||_\infty=|s^\perp|||\chi'||_\infty$, hence
\be ||s_p|^2 div_{\pi_P(S)}g(x)\chi_m(|y|)+(g\cdot s_p)(x)\p_{s^\perp}\chi_m(|y|)|\le C(g,\chi)|s_p|,\ee
where $C(g,\chi)=||g||_\infty||\chi'||_\infty+||\chi_m||_\infty||divg||_\infty$.
Therefore
\be \begin{split}
\left|\int_{G_1(C_m\bs C_{m-1})}\{|s_p|^2 div_{\pi_P(S)}g(x)\chi_m(|y|)+(g\cdot s_p)(x)\p_{s^\perp}\chi_m(|y|)\}dV(z,S)\right|\\
\le C(g,\chi)\int_{G_1(C_m\bs C_{m-1})}|s_p|dV(z,S).\end{split}\ee

Recall that $C_m=B\times B_{\R^{n-d}}(0,m)$, where $B\supset supp\ g$ is compact. By definition of ${\pi_P}_\s(V)$: 
\be {\pi_P}_\s(V)(G_1(B))=\int_{G_1(B\times\R^{n-d})}|\pi_P\circ\pi_S|dV(z,S)=\int_{G_1(B\times\R^{n-d})}|s_p|dV(z,S).\ee
 By hypothesis, ${\pi_P}_\s(V)$ is locally finite, hence $\int_{G_1(B\times\R^{n-d})}|s_p|dV(z,S)<\infty$. As a result, $\int_{G_1(B\times\R^{n-d}\bs C_m)}|s_p|dV(z,S)$ tends to 0 as $m\to \infty$, in particular, $\int_{G_1(C_m\bs C_{m-1})}|s_p|dV(z,S)\to 0$ as $m\to \infty$. Thus by (2.17), the first term of (2.15) tends to 0 as $m\to \infty$. As a result, the second term of (2.15) tends to 0 as well. That is, 
 \be\int_{G_1(C_{m-1})}|s_p|^2 div_{\pi_P(S)}g(x)dV(z,S)\to 0\ as\ m\to\infty. \ee
 
On the other hand, for the second term of (2.15), we have $\sup_{S\in G_1(P)}div_{\pi_P(S)}g(x)\le||D g||_\infty$, hence 
\be \begin{split}|\int_{G_1(C_m)}|s_p|^2 div_{\pi_P(S)}g(x)dV(z,S)|&\le ||Dg||_\infty\int_{G_1(C_m)}|s_p|^2dV(z,S)\\
&\le ||Dg||_\infty\int_{G_1(B\times \R^{n-d})}|s_p|dV(z,S)<\infty,\end{split}\ee
by dominated convergence theorem, and (2.19)
\be 0=\lim_{m\to\infty}\int_{G_1(C_m)}|s_p|^2 div_{\pi_P(S)}g(x)dV(z,S)=\int_{G_1(B\times \R^{n-d})}|s_p|^2 div_{\pi_P(S)}g(x)dV(z,S).\ee
Recall that $|s_p|^2=|\pi_P\circ \pi_S|$, hence
\be 0=\int_{G_1(B\times \R^{n-d})}|\pi_P\circ \pi_S|^2 div_{\pi_P(S)}g(x)dV(z,S)=\int_{G_1(B)}div_Sg(x)d{\pi_P}_{**}(V)(x,S).\ee
This proves that ${\pi_P}_{**}(V)$ is a stationary varifold on $P$.\qed

Next, we give the following lemma, which says that tangent varifolds of a weighted projection are exactly the weighted projections of tangent varifolds:

\begin{lem} Let $V$ be a 1-dimensional varifold in $\R^n$, $P\in G(n,k)$ for some $1\le k\le n$, $x\in [V]$, then for any $r>0$, 
\be {\pi_P}_{**}\circ {\eta_{x,r}}_\s (V)={\eta_{\pi_P(x),r}}_\s \circ {\pi_P}_{**}(V).\ee
In particular, if ${\pi_P}_{**}(V)$ is locally finite near $\pi_P(x)$, then
\be {\pi_P}_{**}(Var\ Tan\ (V,x)) =Var\ Tan\ ({\pi_P}_{**}V,\pi_P(x)).\ee
\end{lem}

\nd This is just by 
definition. Let $V$, $P$, $x\in [V]$ be as in the statement. We have, for any $g\in C_c^\infty(G_1(P))$,
\be \begin{split}\int_{G_1(P)} g(y,S)&d{\pi_P}_{**}\circ {\eta_{x,r}}_\s (V)(y,S)=\int_{G(1,n)}g(\pi_P(z),\pi_P(T))|\pi_P\circ\pi_T|^2 d {\eta_{x,r}}_\s (V)(z,T)\\
&=\int_{G(1,n)}g(\pi_P\circ\eta_{x,r}(z),\pi_P\circ D\eta_{x,r}(T))|\pi_P\circ\pi_{D\eta_{x,r}(T)}|^2|D\eta_{x,r}\circ\pi_T|dV(z,T)\\
&=\frac1r\int_{G(1,n)}g(\pi_P\circ \eta_{x,r}(z),\pi_P(T))|\pi_P\circ\pi_T|^2 dV(z,T)\\
&=\frac1r\int_{G(1,n)}g(\eta_{\pi_P(x),r}\circ\pi_P(z),\pi_P(T))|\pi_P\circ\pi_T|^2 dV(z,T)\\
&=\frac1r\int_{G_1(P)}g(\eta_{\pi_P(x),r}(y),S)d{\pi_P}_{**}(V)(y,S)\\
&=\int_{G_1(P)}g(\eta_{\pi_P(x),r}(y),D\eta_{\pi_P(x),r}(S))|D\eta_{\pi_P(x),r}\circ \pi_S|d{\pi_P}_{**}(V)(y,S)\\
&=\int_{G_1(P)}g(y,S)d{\eta_{\pi_P(x),r}}_\s\circ {\pi_P}_{**}(V)(y,S).
\end{split}\ee
And (2.25) follows by definition of tangent varifolds. \qed

\section{Restriction to stationary 1-varifold with finite projection images}

We want to apply Proposition 2.6 to every stationary 1-varifold and its projections to almost all subspaces $P$. But the condition of locally finite projection image might not be satisfied. For example,

\begin{ex}let $\{L_k\}_{k\in \N}$ be a dense subset of $G(n,1)$, and $\{z_k\}_{k\in \N}=\Z^n\bs\{0\}$. Let $V=\sum V_{L_k+z_k}$. Obviously it is a stationary 1-varifold (it is not hard to verify that $L_k+z_k$ cannot cluster around any point, hence $V$ is locally finite in $\R^n$). However, since the image of $\{z_k\}_{k\in \N}$ under the radial projection to the unit sphere is dense, the image of $V$ under radial projection to the unit sphere $S^{n-1}$ is locally infinite everywhere, which means that for any direction $\sigma\in S^{n-1}$, the image of the projection of $V$ along $\sigma$ to the $n-1$ plane $\sigma^\perp$ is locally infinite around the point 0. 
\end{ex}

However, the locally finite projection is somehow a global property, because for the part of a varifold in any bounded set, it's projection is always finite (since the varifold itself is locally finite). Hence the finiteness of projection images only concerns the behavior of the far away part of a varifold.  The tangent behaviour of a varifold, on the contrary, is a local property and intuitively has nothing to do with what the varifold looks like at infinity. 

As in the above example, if we want to study tangent behavior at the point 0, it is enough to study those parts that pass by 0. There are lines far away, and they results in infinite projection images. But we can simply forget them--if we take away all those lines outside $B(0,1)$, the rest part is still a stationary varifold. This is the purpose of the following proposition. 

\begin{pro}Let $V$ be a 1-dimensional stationary varifold in $\R^n$, then for any $y\in [V]$, there exists a 1-dimensional stationary varifold $W$ in $\R^n$, such that Var Tan$(V,y)=$Var Tan$(W,y)$, and ${\pi_P}_\s(W)$ is locally finite (and hence $\pi_{**}(W)$ is a stationary 1-varifold in $P$) for all linear subspaces (of arbitrary dimension) $P$ of $\R^n$.
\end{pro}

In order to prove the proposition, we need the following lemma:

\begin{lem}Let $V$ be a 1-stationary varifold in $\R^n$, and $y\in [V]$. For any $r>0$, let $V_r$ denote the varifold $V|_{G_1(B(y,r))}$. Then for almost every $r>0$, $|\d V_r|$ is finite and supported on $\partial B(y,r)$, and there exists a $|\d V_r|$ measurable function $\omega_r:\R^n\to S^{n-1}$, such that
\be \d V_r(g)=\int_{\R^n}<g(x)\cdot \omega_r(x)>d|\d V_r|(x)\ee
for any $g\in C_C^\infty(\R^n,\R^n)$. Moreover, 
\be <\omega_r(x),\frac{x}{|x|}>\ge 0\mbox{ for }|\d V_r|-a.e.\ x\in \R^n,\ee
 that is, $\omega_r(x)$ points ourwards from the ball $B(y,r)$.
\end{lem}


\nd


Let $f(x)=|x-y|$. Then $f:\R^n\to \R$ is smooth on $\R^n\bs \{y\}$. We apply Theorem 4.10 (2) of \cite{All72} to our varifold $V$ and the map $f$, and get that, for all most all $r>0$, $|\d V_r|$ is a Radon measure. Then for these $r$, we apply the Riesz representation theorem (cf. \cite{Fe} 2.5.12), and get the existence of a $|\d V_r|$ measurable function $\omega_r:\R^n\to S^{n-1}$, such that
\be \d V_r(g)=\int_{\R^n}<g(x)\cdot \omega_r(x)>d|\d V_r|(x)\ee
for any $g\in C_C^\infty(\R^n,\R^n)$.

To prove that $|\d V_r|$ is supported on $\partial B(y,r)$, just notice that for any vector field $g\in C_C^\infty(\R^n,\R^n)$ with support in $(\partial B(y,r))^C$, 
\be \d V_r(g)=\left\{\begin{array}{cc}\d V(g)&\mbox{ if spt }g\subset B(y,r)\\
0  &\mbox{ if spt }g\subset \overline{B(y,r)}^C\end{array}\right\}=0\ee
since $V$ is stationary.

Since $|\d V_r|$ is locally finite, and supported on a compact set, it is finite.
 
To determine the orientation of the vector field $\omega_r$, note that $V$ is stationary, we can thus apply \cite{All72}  Lemma 4.10 (1) to $V$ and $f(x)=|x|$ and get
\be\begin{split} 0&=\d (V|_{G_1(B(0,r)^C)})(g)+\lim_{h\to 0} \frac 1h\int_{(x,S): |x|\in (r,r+h]} \pi_S(g(x))\cdot \frac{x}{|x|} dV(x,S)\\
&=-\d V_r(g)+\lim_{h\to 0} \frac 1h\int_{(x,S): |x|\in (r,r+h]} <\pi_S(\frac{x}{|x|}), g(x)>dV(x,S),\end{split}\ee
which yields
\be \d V_r(g)=\lim_{h\to 0} \frac 1h\int_{(x,S): |x|\in (r,r+h]}<\pi_S(\frac{x}{|x|}), g(x)>dV(x,S)\ee
for every $g\in C_c^\infty(\R^n,\R^n)$.

Here the vectors $\pi_S(\frac{x}{|x|})$ satisfy that $<\pi_S(\frac{x}{|x|}),\frac{x}{|x|}>\ge 0$ for all $x$ and $S$. We compare with (3.6), and get that for $|\d V_r|$-a.e. $x\in \R^n$, $<\omega_r(x),\frac{x}{|x|}>\ge 0$.\qed

\noindent\textbf{Proof of Proposition 3.2.}

Let $V$ and $y\in [V]$ as in Proposition 3.2. By Lemma 3.3, take an $r \in (0,1]$ such that $|\d V_r|$ is finite and supported on $\partial B(y,r)$, and (3.4) holds. Fix this $r$. For all $x\in \partial B(y,r)$, Set $E_x=\{x+t\omega_r(x),t\ge 0\}$. That is, $E_x$ is the half line issued from the point $x$ with direction $\omega_r(x)$. By (3.5), $E_x\cap B(y,r)=\emptyset$.

Set $V'=\int_{\partial B(y,r)}V_{E_x}d|\d V_r|(x)$. Then it is a 1-varifold, and a simple calculus gives
\be \d V'=-\int_{\partial B(y,r)}<\omega_r(x),\cdot>d|\d V_r|(x)=-\d V_r.\ee

So we set $W=V'+V_r$. Then it is a stationary 1-varifold.

Let $P\in G(n,k)$ for some $k<n$. Then ${\pi_P}_\s (W)={\pi_P}_\s (V_r)+{\pi_P}_\s(V')$. ${\pi_P}_\s(V_r)$ is locally finite because $V_r$ is finite (Remark 2.3); For $V'$, note that $V'=\int_{\partial B(y,r)}V_{E_x}d|\d V_r|(x)$, hence
\be {\pi_P}_\s(V')=\int_{\partial B(y,r)}{\pi_P}_\s(V_{E_x})d|\d V_r|(x)=\int_{\partial B(y,r)}(V_{\pi_P(E_x)})d|\d V_r|(x).\ee

But $|\d V_r|$ is finite, hence ${\pi_P}_\s(V')$ is locally finite. As a result, ${\pi_P}_\s(W)$ is locally finite for all $P\in G(n,k)$, $\forall k<n$.

The last thing for us to verify is that Var Tan$(V,y)=$Var Tan$(W,y)$. So let $f$ be any continuous function with compact support on $G_1(\R^n)$. Since its support is compact, there exists $R>0$ such that $supp\ f\subset G_1(B(0,R))$. Then for any $t<\frac rR$, 
\be\begin{split}
\int f(x,S)d{\eta_{y,t}}_\s V(x,S)&=\int_{x\in B(0,R)}d{\eta_{y,t}}_\s V(x,S)=\frac 1t\int_{\eta_{y,t}(x)\in B(0,R)} f(\eta_{y,t}(x),S)d V(x, S)\\
&=\frac1t\int_{x\in B(y,tR)}f(\eta_{y,t}(x),S)d V(x, S)=\frac1t\int_{x\in B(y,r)}f(\eta_{y,t}(x),S)d V(x, S)\\
&=\frac1t\int_{x\in B(y,r)}f(\eta_{y,t}(x),S)d W(x, S)=\int f(x,S)d{\eta_{y,t}}_\s W(x,S)
\end{split}\ee
since $V$ and $W$ coincide in $B(y,r)$. Hence Var Tan$(V,y)=$Var Tan$(W,y)$ by Definition 1.18 of tangent varifold. This completes the proof of Proposition 3.2. \qed

After the proposition, when we study tangent behavior for stationary 1-varifolds, it is enough to study those with locally finite projection images for almost all directions, and use Proposition 2.6 freely.

\section{Uniqueness of tangent varifolds for stationary 1-varifolds on rectifiable points}

In this section we are going to prove the uniqueness tangent behaviour for stationary 1-varifolds in $\R^n$ at a particular class of points: the rectifiable points, which are defined as follow:

\begin{defn}Let $V$ be a 1-varifold in $\R^n$. A point $x\in [V]$ is called a rectifiable point, if Var Tan$(V,x_0)$ contains at least one rectifiable element.
\end{defn}

In the rest of this section, we only consider projections on $n-1$-dimensional affine planes. Note that if $P$ is parallel to $P'$, then the projection images or weighted projections of any $1-$varifold on them are the same. Hence we do not distinguish them. And for any $n-1$-dimensional plane $P$, let $\sigma\in S^{n-1}$ be the unit vector orthogonal to $P$. Then we also call projections on $P$ the projection to the direction $\sigma$.

\begin{thm}[Uniqueness of tangent varifold at rectifiable points] Let $V$ be a stationary 1-varifold in $\R^n$. Then for any rectifiable point $x_0\in [V]$, Var Tan$(V,x_0)$ contains only one element.
\end{thm}

\begin{rem}Theorem 4.2 can be deduced directly from Theorem 5.1 for general case in the next section. We prove it separately, because here for rectifiable case the proof is more geometric; and as we will see in Lemma 4.21, to ensure that two conic rectifiable 1-varifolds are the same, we only need that they have the same weighted projections on $n-1$-subspaces in a set of positive measure in the Grassmannian $G(n,n-1)$. But in the proof for non-rectifable case in Section 5, we need that they have same weighted projections in almost all $2$-subspaces (and hence in almost all $n-1$-subspaces).
\end{rem}

\noindent\textbf{Proof of Theorem 4.2.} Thanks to Proposition 3.2, it is enough to prove it for all stationary 1-varifolds in $\R^n$ with locally finite projection images to all directions.

We are going to prove by recurrence on ambient dimensions. As we will see, the recurrence argument can only apply when $n\ge 2$ (cf. Remark 4.61). So let us first check it for $n=1,2$:

For $n=1$: by \cite{Sim83} Page 243 we know that any element in Var Tan$(V,x_0)$ for a 1-dimensional stationary varifold is conical. As a result, Var Tan$(V,x_0)=\{\theta(x_0)V_\R\}$, where $\theta(x_0)$ is the density of $V$ at $x_0$.;

For $n=2$: this is due to Hirsch \cite{HP09}, where he proves that any 1-dimensional stationary varifold (not necessarily rectifiable) admits unique tangent varifold everywhere.

Now suppose that Theorem 4.2 is true for some $n\in \N$ and $n\ge 2$, and we will prove it for $n+1$:

Let $V$ be any stationary 1-rectifiable varifold in $\R^{n+1}$ such that ${\pi_P}_\s(V)$ is locally finite for all $P\in G(n+1,n)$, and let $x_0\in[V]$ be a rectifiable point. Let $W_0\in$ Var Tan$(V,x_0)$ be rectifiable. By Lemma 2.23, for all $P\in G(n+1,n)$, ${\pi_P}_{**}(W_0)$ is a rectifiable element in Var Tan$({\pi_P}_{**}(V), \pi_P(x_0))$. Since ${\pi_P}_\s(V)$ is locally finite for all $P\in G(n,n+1)$, by Proposition 2.6, ${\pi_P}_{**}(V)$ is a stationary 1-varifold in $P$, which is of dimension $n$. Hence by hypothesis of recurrence, Var Tan$({\pi_P}_{**}(V))$ contains only one element. This yields
\be {\pi_P}_{**}(W_0)={\pi_P}_{**}(W_1),\forall P\in G(n+1,n),\forall W_1\in \mbox{Var Tan}(V,x_0).\ee
In other words, any element in Var Tan$(V,x_0)$ has the same weighted projection as $W_0$ for all directions.


Due to the above projection property, and the fact that tangent varifolds for 1-stationary varifolds are conical (cf. \cite{Sim83} P243), let us first study the projection of conical 1-varifolds.
%
%

%

So let $W$ be any conical 1-varifolds. Then there is a Radon measure $\mu$ on $S^n$, such that 
\be W=\int_{z\in S^n}V_{R_z}d\mu(z),\ee
Where for any point $z\in S^n$, $R_z$ is the half line issued from the origin and passing through the point $z$.

Let $\phi(z)=\mu(\{z\})$. Note that for any $z\in S^n$, $\theta_1(W,z)=\phi(z)$.

\begin{rem}For a general conical 1-varifold $W$, we have $\mu\ge \sum_{z\in S^n}\phi(z)\d_z$. In particular, if $W$ is a 1-rectifiable varifold, then $\phi:S^n\to [0,\infty)$ is a map which is zero except for a countable subset of $S^n$, and $\mu=\sum_{z\in S^n}\phi(z)\d_z$, $W=\sum_{z\in S^n}\phi(z)V_{R_z}$. 
\end{rem}

\begin{lem}Let $W=\int_{z\in S^n}V_{R_z}d\mu(z)$ be a conical 1 varifold in $\R^{n+1}$, where $\mu$ is a Radon measure on $S^n$. Let $\phi(z)=\mu(\{z\})$ for $z\in S^n$. Let $y\in S^n$. Set $W_r=W|_{G_1(B(y,r))}$ for any $r>0$. Then

$1^\circ$, for all $r<\frac 12$, and all $P\in G(n+1,n)$,
\be {\theta_1}_*({\pi_P}_{**}W_r,\pi_P(y))\ge |\pi_P(y)|\phi(y);\ee

$2^\circ$  Let $\e>0$ be small, and let $r_0\in(0,10^{-5})$ be such that $\mu(S^n\cap B(y,r))<\phi(y)+\e$. Then for any $r<r_0$, $P\in G(n,n+1)$ with $|\pi_P(y)|>10 r$, we have
\be {\theta_1}^*({\pi_P}_{**}W_r,\pi_P(y))\le \frac32\e+|\pi_P(y)|\phi(y).\ee
\end{lem}

\nd 

Denote by $a=\phi(y)=\theta_1(W,y)\ge 0$ for short. 

For $1^\circ$, it is trivial if $\phi(z)=0$ or $|\pi_P(y)|=0$, so suppose that $a=\phi(y)>0$, and $|\pi_P(y)|\ne 0$. We have, for any $t\in (0,|\pi_P(y)|r)$, 
\be\begin{split} 
&|{\pi_P}_{**}(aV_{[(1-r)y,(1+r)y]})|(B_P(\pi_P(y),t))\\
=&a\int_{\{(x,S): \pi_P(x)\in B_P(\pi_P(y),t)\}}|\pi_P\circ \pi_S|^2dV_{[(1-r)y,1+r)y]}\\
=&a\int_{\{x\in [(1-r)y,1+r)y],|\pi_P(x-y)|<t\}}|\pi_P(y)|^2d\H^1(x)\\
=&a\int_{\{x\in [(1-\frac{t}{|\pi_P(y)|})y,1+\frac{t}{|\pi_P(y)|})y]\}}|\pi_P(y)|^2d\H^1(x)=2ta|\pi_P(y)|.
\end{split}\ee

This is true for all $t$ small, hence 
\be \begin{split}{\theta_1}({\pi_P}_{**}(aV_{[(1-r)y,(1+r)y]}), \pi_P(y))&=\lim_{t\to 0} |{\pi_P}_{**}(aV_{[(1-r)y,(1+r)y]})|(B_P(\pi_P(y),t))/2t\\
&=a|\pi_P(y)|=|\pi_P(y)|\phi(y),\end{split}\ee
and therefore
\be {\theta_1}_*({\pi_P}_{**}(W_r),\pi(P(y)))\ge {\theta_1}({\pi_P}_{**}(aV_{[(1-r)y,(1+r)y]}), \pi_P(y))=|\pi_P(y)|\phi(y),\ee
which yields $1^\circ$.

Now let us look at $2^\circ$. Let $P\in G(n+1,n)$. For any $t>0$ small enough, we have
\be\begin{split}
& |{\pi_P}_{**}(W_r)|B_P(\pi_P(y), t)=|{\pi_P}_{**}(aV_{[(1-r)y,(1+r)y]})|(B_P(\pi_P(y),t))+\\
&\ \ \ \ \ \ \int_{z\in B(y,r)\cap S^n\bs \{y\}}d\mu(z)\int_{\{(x,S):\pi_P(x)\in B_P(\pi_P(y),t)\}}|\pi_P\circ \pi_S|^2 dV_{R_z\cap B(y,r)}\\
 &=2ta|\pi_P(y)|+\int_{z\in B(y,r)\cap S^n\bs \{y\}}d\mu(z)\int_{\{(x,S):\pi_P(x)\in B_P(\pi_P(y),t)\}}|\pi_P\circ \pi_S|^2 dV_{R_z\cap B(y,r)}\\
 &\le 2ta|\pi_P(y)|+\int_{z\in B(y,r)\cap S^n\bs \{y\}}d\mu(z)\int_{\{(x,S):\pi_P(x)\in B_P(\pi_P(y),t)\}}|\pi_P\circ \pi_S|^2 dV_{[(1-2r)z,(1+2r)z]}\\
 &= 2ta|\pi_P(y)|+\int_{z\in B(y,r)\cap S^n\bs \{y\}}d\mu(z)\int_{\{x\in[(1-2r)z,(1+2r)z], |\pi_P(x-y)|<t\}}|\pi_P(z)|^2 d\H^1(x).\end{split}\ee
 
 Note that $z\in B(y,r)\cap S^n$, hence $|\pi_P(y)|-r\le |\pi_P(z)|\le |\pi_P(y)|+r$. Also, the length of the segment 
 \be\H^1([(1-2r)z,(1+2r)z]\cap \{x: |\pi_P(x-y)|<t\})\le \frac{2t}{|\pi_P(z)|}\le \frac{2t}{|\pi_P(y)|-r},\ee
 hence
 \be\int_{\{x\in[(1-2r)z,(1+2r)z], |\pi_P(x-y)|<t\}}|\pi_P(z)|^2 d\H^1(x)\le (|\pi_P(y)|+r)^2 \frac{2t}{|\pi_P(y)|-r},\ee
 and therefore
 \be |{\pi_P}_{**}(W_r)|B_P(\pi_P(y), t)\le 2ta|\pi_P(y)|+\int_{z\in B(y,r)\cap S^n\bs \{y\}}d\mu(z)(|\pi_P(y)|+r)^2 \frac{2t}{|\pi_P(y)|-r}.\ee
 
 Now if $|\pi_P(y)|\ge 10 r$, then
 \be (|\pi_P(y)|+r)^2 \frac{2t}{|\pi_P(y)|-r}\le 3t|\pi_P(y)|,\ee
 and hence
 \be\begin{split} |{\pi_P}_{**}(W_r)|B_P(\pi_P(y), t)&\le 2ta|\pi_P(y)|+3t(\mu(B(y,r)\cap S^n\bs \{y\}))\\
 \le 2ta|\pi_P(y)|+3t\e,
\end{split}\ee
which yields
\be {\theta_1}^*({\pi_P}_{**}W_r,\pi_P(y))\le \frac32\e+ |\pi_P(y)|\phi(y),\forall r<r_0.\ee
\qed

Now we let us go back to $W_0$ and $W_1$ in (4.4), where $W_0$ is a rectifiable element in Var Tan$(V,x_0)$, and $W_1$ is any element in Var Tan$(V,x_0)$. Let $\mu_i$ be the corresponding Radon measures of $W_i$ on $S^n$, and $\phi_i:S^n\to [0,\infty)$, $\phi_i(z)=\mu_i(\{z\})$ for $z\in S^n$, $i=0,1$. Then 
\be W_i=\int_{z\in S^n}V_{R_z}d\mu_i(z), i=1,2.\ee

In particular, since $W_0$ is rectifiable, by Remark 4.6,
$\mu_0=\sum_{z\in S^n}\phi_0(z)\d_z$, and $W_0=\sum_{z\in S^n}\phi_0(z)V_{R_z}$. (Note that we do not have this for $W_1$ yet, since we do not know whether it is rectifiable.)

We want to prove that $\mu_0=\mu_1$. Let us first prove that $\mu_0\le \mu_1$, i.e. $\phi_0\le \phi_1$ for all $z\in S^n$. To see this, it is enough to consider those $z\in S^n$ with $\phi_0(z)>0$.

\begin{lem}$\phi_0\le \phi_1$ on $S^n$.
\end{lem}

We prove by contradiction. Suppose there is a $y\in S^n$ such that $\phi_0(y)>\phi_1(y)$. Then $\phi_0(y)>0$, and there exists $\e\in (0,10^{-5})$ such that $\phi_1(y)<\phi_0(y)-5\e$.

 By continuity of Radon measure, there is an $r\in (0,10^{-5})$ such that $\mu_1(B(y,r))<\phi_1(y)+\e$. Fix this $r$. Denote by $W^i_r=W_i|_{G_1(B(y,r))}$, $i=0,1$. Then by Lemma 4.7, we have
 \be \theta_1({\pi_P}_{**}W_0,\pi_P(y))\ge {\theta_1}_*({\pi_P}_{**}W^0_r,\pi_P(y))\ge |\pi_P(y)|\phi_0(y),\ee
 and for $P\in G(n+1,n)$ with $|\pi_P(y)|>10^{-4}>10r $, we have
 \be {\theta_1}^*({\pi_P}_{**}W^1_r,\pi_P(y))\le \frac 32\e+|\pi_P(y)|\phi_1(y).\ee
 
 Thus, for $P\in G(n+1,n)$ with $|\pi_P(y)|>\frac 12$, 
 \be {\theta_1}^*({\pi_P}_{**}W^1_r,\pi_P(y))\le 3\e |\pi_P(y)|+|\pi_P(y)|\phi_1(y)=(3\e+\phi_1(y))|\pi_P(y)|< (\phi_0(y)-\e)|\pi_P(y)|.\ee
 But $W_0$ and $W_1$ have the same weighted projections, hence by (4.22),
 \be \theta_1(({\pi_P}_{**}W_1,\pi_P(y))=\theta_1({\pi_P}_{**}W_0,\pi_P(y))\ge |\pi_P(y)|\phi_0(y).\ee
 
So
\be \theta_1(({\pi_P}_{**}W_1,\pi_P(y))-{\theta_1}^*({\pi_P}_{**}W^1_r,\pi_P(y))\ge \e|\pi_P(y)|\ge \frac{\e}2>0\ee
for all $P\in G(n+1,n)$ with $|\pi_P(y)|>\frac12.$

Set $W'=W_1|_{G_1(B(y,r)^C)}$, then (4.26) yields
\be {\theta_1}_*({\pi_P}_{**}(W'), \pi_P(y))\ge \frac\e2\ee
for all $P\in G(n+1,n)$ with $|\pi_P(y)|>\frac 12$. As a result,
\be {\theta_1}_*({\pi_P}_\s (W'), \pi_P(y))\ge \frac\e2\ee
for all $P\in G(n+1,n)$ with $|\pi_P(y)|>\frac 12$, since by definition,
$|{\pi_P}_\s (W')|\ge |{\pi_P}_{**}(W')|.$

Note that $W'$ is a 1-varifold with support outside $B(y,r)$. 

Denote by $f$ the shortest distance projection from $B(y,r)^C\to \partial B(y,r)$, that is, $f(x)=y+r\frac{x-y}{|x-y|}$ for $x\in B(y,r)^C$. 


We claim that, for any $P\in G(n+1,n)$ with $|\pi_P(y)|\ge \frac 12$, 
\be \theta_{1*}(f_\s(W'),y+r P^\perp)+\theta_{1*}(f_\s(W'),y-r P^\perp)\ge \frac14 \theta_{1*}({\pi_P}_\s(W'),\pi_P(y)),\ee
where $P^\perp$ is a unit normal vector to $P$.
 
To prove Claim 4.29, we fix any $P\in G(n+1,n)$ such that $\pi_P(y)\ne 0$, and let $e\in \R^{n+1}$ be a unit normal vector to $P$. Set $Q+=\{x: <x-y,e>\ge 0\}$ and $Q_-=\{x:<x-y,e>\le 0\}$. Also fix any $t<10^{-5}$. We want to estimate
$ |{\pi_P}_\s(W'_{G_1(Q_+)})|(B(\pi_p(y),tr))$ and $|f_\s(W')|(B(y+re,tr))$. Since $W'$ coincides with a conic varifold outside $B(y,r)^C$, it is enough to restrict ourselves to studying the case when $W'=V_{R_z}|_{B_1(G(y,r)^C)}$, the restriction to $B(y,r)^C$ of a varifold generated by the ray $R_z$ for some $z\in \R^n\bs \{0\}$. In this case, by Remark 1.16, 
\be |{\pi_P}_\s(W'_{G_1(Q_+)})|(B(\pi_p(y),tr))=\H^1(\pi_P(R_z\cap Q_+\bs B(y,r))\cap B(\pi_p(y),tr))=\H^1(\pi_P(R_z\cap D_t\bs B(y,r))),\ee 
where $D_t$ denotes the half cylinder ${\pi_P}^{-1}(B(\pi_p(y),tr))\cap Q_+$, and 
\be\begin{split} |f_\s(W')|(B(y+re,tr))&=\H^1(f(R_z\bs B(y,r))\cap B(y+re,tr))\\
&=\H^1(f(R_z\cap f^{-1}(B(y+re,tr)))).\end{split}\ee
Set $C_t:=\{x: <x-y,e>\ge \sqrt{1-t^2} |x-y|\}=\{x\in Q_+:|\pi_P(x-y)|\le t|x-y|\}$, which is a cone centered at $y$, then $C_t\cap\partial B(y,r)=D_t\cap \partial B(y,r)$, and $C_t\bs B(y,r)\subset f^{-1}(B(y+re,2tr))$. Hence
\be |f_\s(W')|(B(y+re,2tr))\ge \H^1(f(R_z\cap C_t\bs B(y,r))).\ee

Also note that, since $|\pi_P(y)|\ge \frac 12$, $0\not\in C_t$ and $0\not\in D_t$. So since $R_z$ is a half line, and its endpoint (the origin) does not belong to either $D_t$ or $C_t$, the intersections $R_z\cap D_t\bs B(y,r)$, and $R_z\cap C_t\bs B(y,r)$, are either segments, or unions of two disjoint segments, with endpoints on $R_z$.

Set $g:\R^{n+1}\to \R^{n+1}$, $g(x)= y+\frac{r}{|<x-y,e>|}(x-y)$. For any segment $[z_0,z_1]\subset B(y,r)^C$ such that $f([z_0,z_1])\subset B(\pi_p(y),tr))$, $[z_0,z_1]$, and the images $f[z_0,z_1]$ and $g[z_0,z_1]$, are contained in the cone $C_t$ and the affine plane $\Pi$ passing by $y,z_0$ and $z_1$. Moreover, the image $f([z_0,z_1])$ is the arc of great circle with endpoints $f(z_0)$ and $f(z_1)$ on the sphere $\partial B(y,r)$, and $g[z_0,z_1]$ is a segment parallel to $P$ with distance $r$. So if we look at everything in the plane $\Pi$, $f([z_0,z_1])$ is the arc of the circle $\partial B(y,r)\cap \Pi$, and $g[z_0,z_1]$ is a segment parallel to $P$ with distance at most $\sqrt{r^2+(tr)^2}$ to $P\cap \Pi$.

Hence 
\be\begin{split}
\frac{\H^1(g[z_0,z_1])}{\H^1(f[z_0,z_1])}&\le\sup_{\theta_1,\theta_2\in [-\arcsin t,\arcsin t]}\frac{(\tan\theta_1-\tan\theta_2)\sqrt{r^2+(tr)^2}}{(\theta_1-\theta_2)r}\\
&\le\sqrt{r^2+(tr)^2}\sup_{\theta\in [-\arcsin t,\arcsin t]}\tan '(\theta)\le 2.
\end{split}\ee

Therefore for any $t<10^{-5}$, and for any segment $[z_0,z_1]\subset C_t$, 
\be \H^1(g[z_0,z_1])\le 2\H^1(f[z_0,z_1]).\ee

Since $R_z\cap C_t\bs B(y,r)$ is either a segment, or a union of two disjoint segments, by (4.34), 
\be \H^1(g(R_z\cap C_t\bs B(y,r)))\le 2 \H^1(f(R_z\cap C_t\bs B(y,r))).\ee

On the other hand, let us look at $R_z\cap D_t\bs B(y,r)$. By the argument above, it is also a segment or a union of two segments. 

If $[R_z\cap D_t]\cap B(y,r)=\emptyset$, then $R_z\cap D_t\bs B(y,r)=R_z\cap D_t$ is a segment. Denote $R_z\cap D_t$ by $[z_0,z_1]$, and denote by $z'$ its midpoint. Then since $r$ is small, and $|\pi_P(y)|\ge \frac 12$, the endpoint of $R_z$ (the origin) is outside $D_t$, hence $\pi_P(z_0),\pi_P(z_1)\in \partial B(\pi_P(y),rt)$. As a result, $\pi_P([z_0,z_1])$ is a chord of the ball $B(\pi_P(y),rt)$. Therefore the two segments $[\pi_P(y),\pi_P(z')]$ and $[\pi_P(z_0),\pi_P(z_1)]$ are mutually perpendicular, unless $\pi_P(z')=\pi_P(y)$.

In case $\pi_P(z')\ne\pi_P(y)$, let $e_1=\frac{\overrightarrow{\pi_P(z_0)\pi_P(z_1)}}{|\pi_P(z_0)-\pi_P(z_1)|}\in P$, and $e_2=\frac{\overrightarrow{\pi_P(y)\pi_P(z')}}{|\pi_P(y)-\pi_P(z')|}\in P$. Then $\{e,e_1,e_2\}$ forms an orthonormal set, and $y,e,R_z,\pi_P(R_z),f(R_z),g(R_z)$ are all contained in the 3-dimensional subspace $E$ generated by these three vectors. So we will work only in this subspace with coordinates under the orthonormal basis $\{e,e_1,e_2\}$.

As a result, there exists $b,c,a_0,a_1\ge 0$ with $c^2+b^2=(rt)^2$, such that
\be \pi_P(z')=\pi_P(y)+(0,0,b), \pi_P(z_0)=\pi_P(y)+(0,-c,b), \pi_P(z_1)=\pi_P(y)+(0,c,b),\ee
and
\be z_0=y+(a_0,-c,b),z_1=y+(a_1,c,b),z'=y+(\frac{a_0+a_1}{2},0,b).\ee

Immediately, we get  
\be\H^1(\pi_P(R_z\cap D_t))=|\pi_P(z_1)-\pi_P(z_0)|=2c.\ee

Note that when $t$ is small, $a_0$ and $a_1$ are both positive. Otherwise, since $\frac{a_0+a_1}{2}>0$, $a_0$ and $a_1$ are of different signs. But $a_i=<z_i-y,e_0>$, for $i=0,1$, hence there exists a point $z_2\in [z_0,z_1]$ such that $<z_2-y,e_0>=0$, and thus $|z_2-y|=|\pi_P(z_2-y)|\le rt$, hence $z_2\in B(y,r)$, and this contradicts the fact that $[z_0,z_1]\cap B(y,r)=\emptyset$.

Let us now look at $C_t\cap R_z$. We want to estimate $g(C_t\cap R_z\bs B(y,r))$. Recall that we are in the case where $R_z\cap D_t$ does not meet $B(y,r)$. Note that $C_t\bs D_t\subset B(y,r)^C$, hence $R_z\cap C_t$ does not meet $B(y,r)$ either. Also, since $|\pi_P(y)|\ge \frac 12$, and $t<10^{-5}$, $r<10^{-5}$, thus $0\not\in C_t$. Hence $R_z\cap C_t$ is a segment, with endpoints on $\partial C_t$. Therefore $g(R_z\cap C_t)$ is a chord $I$ of the ball $B_t=(P+y+re)\cap C_t$ with radius $r\tan(\arcsin t)=\frac{rt}{\sqrt{1-t^2}}$.

Since $[z_0,z_1]\subset D_t\bs B(y,r)\subset C_t\bs B(y,r)$, $[g(z_0),g(z_1)]=g([z_0,z_1])\subset I$. We calculate $g(z_0)=y+(r,\frac{-rc}{a_0},\frac{rb}{a_0})$ and $g(z_1)=y+(r,\frac{rc}{a_1},\frac{rb}{a_1})$. Without loss of generality, suppose that $a_1>a_0$. Then
\be\begin{split} 
dist&(y+re,I)\le dist(y+re,[g(z_0),g(z_1)])= dist((0,0), [(\frac{-rc}{a_0},\frac{rb}{a_0}),(\frac{rc}{a_1},\frac{rb}{a_1})])\\
&\le dist((0,0), [(\frac{-rc}{a_0},\frac{rb}{a_0}),(\frac{rc}{a_0},\frac{rb}{a_0})])=\frac{r}{a_0}dist(\pi(y),[\pi(z_0),\pi(z_1)])=\frac{rb}{a_0}.
\end{split}\ee

Note that $z_0\in C_t\bs B(y,r)$, hence $a_0=<z_0-y,e>\ge\sqrt{1-t^2}|z_0-y|\ge \sqrt{1-t^2}r$. Therefore $\frac{r}{a_0}\le \frac{1}{\sqrt{1-t^2}}$. As a result, $I$ is a chord of a ball of radias $\frac{rt}{\sqrt{1-t^2}}$, and its distance to the center of the ball is less or equal to $\frac{b}{\sqrt{1-t^2}}$. Hence the length of the chord $I$ is larger or equal to $2\sqrt{(\frac{rt}{\sqrt{1-t^2}})^2-(\frac{b}{\sqrt{1-t^2}})^2}=\frac{2c}{\sqrt{1-t^2}}$, since $c^2+b^2=(rt)^2$. Hence $\H^1(I)\ge \frac{2c}{\sqrt{1-t^2}}\ge 2c$. But $I=g(R_z\cap C_t)$,hence
\be \H^1(g(R_z\cap C_t)\ge 2c=\H^1(\pi_P(R_z\cap D_t)),\ee
and by (4.34),
\be \H^1(f(R_z\cap C_t)\ge \frac12 \H^1(\pi_P(R_z\cap D_t)).\ee

The case where $\pi_P(z')=\pi_P(y)$ is just a degenerated case, and the proof is the same as above.

Now let us discuss the case where $R_z\cap D_t\cap B(y,r)\ne \emptyset$. In this case,  $R_z\cap D_t\cap B(y,r)^C$ is either a segment $[z_0,z_1]$, or a union of two disjoint segments $[z_0,z_1]\cup [z_0', z_1']$, where $z_0,z_0'\in\partial B(y,r)$ and $z_1,z_1'\in B(y,r)^C\cap \partial D_t$. It is not hard to see that in these two cases, $R_z\cap C_t\bs B(y,r)$ is also either a segment $[z_0,z_2]$, or a union of two disjoint segments $[z_0,z_2]\cup[z_0',z_2']$ respectively, with $z_2,z_2'\in \partial C_t$.

Let $a_0=<z_0-y,e>$. Then $a_0\le |z_0-y|=r$.
%

Let $h: C_t\to C_t$, $h(x)= y+\frac{a_0}{<x-y,e>}(x-y)$. Then $h(z_0)=z_0$ and for all segments in $C_t\bs B(y,r)$, 
\be \H^1(h(I))\le \H^1(g(I))\ee
since $a_0\le r$. Note that for any $x\in C_t$, $<x-y,e>$ is always positive. 

We claim that \be \H^1(\pi_P[z_0,z_1])\le \H^1(h[z_0,z_2]).\ee

To prove the claim, we have 2 cases:

\medskip

\noindent\textbf{Case 1}: $<z_2-y,e> \le a_0$. 

In this case we have 
\be \begin{split}&|h(z_2)-h(z_0)|=\frac{a_0}{<z_2-y,e>}(z_2-y)-(z_0-y)\\
&=\frac{a_0}{<z_2-y,e>}(\pi_P(z_2-y)+<z_2-y,e>e)-(\pi_P(z_0-y)+a_0e)\\
&=\frac{a_0}{<z_2-y,e>}\pi_P(z_2-y)-\pi_P(z_0-y)\\
&=\pi_P(z_2-z_0)+(\frac{a_0}{<z_2-y,e>}-1)\pi_P(z_2-y).
\end{split}\ee

Note that $z_2\in \partial C_t\cap B(y,r)^C$, and $z_0\in D_t$, hence $|\pi_P(z_2-y)|=t |z_2-y|\ge tr>|\pi_P(z_0-y)|$, hence in the triangle $\Delta_{\pi_P(y)\pi_P(z_0)\pi_P(z_2)}$, the edge $[\pi_P(y),\pi_P(z_0)]$ is shorter than the edge $[\pi_P(y),\pi_P(z_2)]$. Therefore the angle $\angle_{\pi_P(y)\pi_P(z_2)\pi_P(z_0)}< \frac\pi 2$, and hence $<\pi_P(z_2)-\pi_P(y),\pi_P(z_2)-\pi_P(z_0)>>0$. As a result, 
\be |\pi_P(z_2-z_0)+(\frac{a_0}{<z_2-y,e>}-1)\pi_P(z_2-y)|\ge |\pi_P(z_2-z_0)|,\ee
since $\frac{a_0}{<z_2-y,e>}-1>0$.

Combine with (4.44), we get
\be |h(z_2)-h(z_0)|\ge |\pi_P(z_2-z_0)|\ge |\pi_P(z_1-z_0)|,\ee
which yields (4.43);

\medskip

\noindent\textbf{Case 2}: $<z_2-y,e> > a_0$.

In this case we also have $<z_1,e>>a_0$.

Note that $z_1\in [z_0,z_2]$, $|\pi_P(z_0-y)|< rt$ and $|\pi_P(z_1-y)|=rt$. Hence $h(z_1)\in (z_0,h(z_2)]\cap (y+a_0e, y+a_0e+\frac{a_0}{r}\pi_P(z_1-y)]$. Now the points $y_0=y+a_0e, h(z_0),h(z_1),h(z_2)$ all belong to the intersection $\overline B_0$ of $C_t$ and the affine plane $y+a_0e+Q$ (unless $y\in \pi_P(R_z)$, which is a degenerate trivial case) where $Q=Vect\{\pi_P(z),\pi_P(z_1-y)\}$ is the linear plane generated by the two vectors $\pi_P(z)$ and $\pi_P(z_1-y)$. Thus $\overline B_0$ is a planar disc, centered on $y_0=y+a_0e$, and perpendicular to $e$. Denote by $r_0=a_0\frac{t}{\sqrt{1-t^2}}$ its radius.

Set $\xi_1=y_0+\pi_P(z_1-y)$. Note that $\xi_1\in y_0+Q$, and $|\xi_1-y_0|=|\pi_P(z_1-y)|=rt$. On the other hand, $z_0\in C_t\bs B(y,r)$, hence $a_0=<z_0,e>\ge \sqrt{1-t^2}r$ and therefore
\be |\xi_1-y_0|=rt\le t\frac{a_0}{\sqrt{1-t^2}}=r_0,\ee
hence $\xi_1\in \overline B_0$. Also note that 
\be |h(z_1)-y_0|=\frac{a_0}{<z_1,e>}|\pi_P(z_1-y)|\le |\pi_P(z_1-y)|=|\xi_1-y_0|.\ee

After the discussion above, we have, $h(z_0)=z_0\in \overline B_0\bs \{y_0\}$, $\xi_1=y+a_0e+\pi_P(z_1-y)\in \overline B_0$, $h(z_1)\in (y_0, \xi_1]$.  Also note that $\H^1(\pi_P[z_0,z_1])=|z_0-\xi_1|$, and $h(z_2)$ is the intersection of $\partial B_0$ and the half-line $R_{z_0h(z_1)}$.

Let $\xi_2$ be the intersection of $\partial B_0$ and the half-line $R_{y_0\xi_1}$. Then we have $\angle_{y_0\xi_2h(z_2)}=\angle_{y_0h(z_2)\xi_2}$

By definition of $z_0$ and $z_1$, we have $|z_0-y_0|=\pi_P(z_0-y)|\le rt=|\pi_P(z_1-y)|=|\xi_1-y_0|$. As a result, in the triangle $\Delta_{z_0\xi_1y_0}$ we have $\angle_{z_0\xi_1y_0}\le \angle_{\xi_1z_0y_0}$. This means that $\angle_{z_0\xi_1y_0}\le \frac\pi 2$, and hence $\angle_{z_0\xi_1\xi_2}\ge \frac\pi 2$. Thus in the triangle $ \Delta_{z_0\xi_1\xi_2}$, we have 
\be |z_0-\xi_2|>|z_0-\xi_1|.\ee

Thus we have $|z_0-\xi_2|\ge |z_0-\xi_1|=\H^1(\pi_P[z_0,z_1])$.

 Since $\xi_1\in[y_0,\xi_2]$, $\angle_{y_0\xi_1h(z_2)}\ge \angle_{y_0\xi_2h(z_2)}$ and $\angle_{y_0h(z_2)\xi_1}\le\angle_{y_0h(z_2)\xi_2}$. Also, $\angle_{y_0h(z_2)\xi_1}\ge \angle_{h(z_1)h(z_2)\xi_1}=\angle_{z_0h(z_2)\xi_1}$ since $h(z_1)\in [y_0,\xi_1]$, and $\angle_{z_0\xi_1h(z_2)}\ge\angle_{h(z_1)\xi_1h(z_2)}=\angle_{y_0\xi_1h(z_2)}$ since $h(z_1)\in[z_0,h(z_2)]$. Altogether, we have
\be \angle_{z_0h(z_2)\xi_1}\le\angle_{z_0\xi_1h(z_2)},\ee
and hence in the triangle $\Delta_{z_0h(z_2)\xi_1}$, $|z_0-\xi_1|\le |z_0-h(z_2)|$, that is, again,
\be \H^1(\pi_P[z_0,z_1])\le \H^1(h[z_0,z_2]).\ee

Hence both Case 1 and Case 2 give us (4.43).

\medskip

Similarly, we have
\be \H^1(\pi_P[z_0',z_1'])\le \H^1(h[z_0',z_2']).\ee

Hence by (4.34) and (4.42),
\be \begin{split}
\H^1(\pi_P(R_z\cap D_t\bs B(y,r))&\le \H^1(h(R_z\cap C_t\bs B(y,r))\\
&\le \H^1(g(R_z\cap C_t\bs B(y,r))\le 2\H^1(f(R_z\cap C_t\bs B(y,r)),\end{split}\ee
which yields again
\be \H^1(f(R_z\cap C_t\bs B(y,r))\ge\frac12 \H^1(\pi_P(R_z\cap D_t\bs B(y,r))).\ee

Combine with (4.41) (the case where $R_z\cap D_t$ and $R_z\cap C_t$ do not meet $B(y,r)$), we get that (4.54) is true in all cases. By (4.30) and (4.32), we have
\be |f_\s(W')|(B(y+re,2tr))\ge \frac12 |{\pi_P}_\s(W'_{G_1(Q_+)})|(B(\pi_P(y),tr)),\ee
which yields
\be \theta_{1*}(|f_\s(W')|,y+re)\ge \frac14 \theta_{1*}(|{\pi_P}_\s(W'_{G_1(Q_+)})|,\pi_P(y)).\ee

Similarly we get
\be \theta_{1*}(|f_\s(W')|,y-re)\ge \frac14 \theta_{1*}(|{\pi_P}_\s(W'_{G_1(Q_-)})|,\pi_P(y)),\ee
 
and altogether we get the claim (4.29).

By (4.28), we get that for $P\in G(n+1,n)$ with $|\pi_P(y)|\ge\frac12$, we have
\be \theta_{1*}(f_\s(W'),y+r P^\perp)+\theta_{1*}(f_\s(W'),y-r P^\perp)\ge \frac\e 8>0.\ee

Thus, Borel measurable set $A=\{\xi\in S^n: {\theta_1}_*(|f_\s(W')|,y+r\xi) >0\}$ is of positive $\H^n$-measure. 

But note that 
\be \begin{split}|f_\s(W')|(S^n)\le |f_\s(W_1)|(S^n)&=\int_{z\in S^n}|f_\s(V_{R_z})|(S^n)d\mu_1(z)\\
&=\int_{z\in S^n}|V_{f(R_z)}|(S^n)d\mu_1(z)=\int_{z\in S^n}\H^1(f(R_z))d\mu_1(z)\\
&\le \int_{z\in S^n}2\pi d\mu_1(z)=2\pi\mu_1(S^n),\end{split}\ee
hence $|f_\s(W')|$ is a finite measure on $S^n$.

But we know that its lower 1-density is positive on a set $A\subset S^n$ of positive $\H^n$ measure, which means that its lower $n$-density is infinite on $A$, since $n\ge 2$. By differentiation theorem between Radon measures, $A$ should of $\H^n$ measure zero, which leads to a contradiction. This finish our proof of Lemma 4.21. \qed

After Lemma 4.21, we know that the measures on $S^n$ satisfy
\be \mu_1\ge \sum_{z\in S^n}\phi_1(z)\d_{z}\ge \sum_{z\in S^n}\phi_0(z)\d_z=\mu_0.\ee

But since both $W_0$ and $W_1$ are tangent varifolds of $V$ on a same point $x_0$, they should admit the same 1-density (which equals $\theta(V,x_0)$) at the origin. Hence $\mu_0(S^n)=\mu_1(S^n)$. Combine with (4.60), we get $\mu_0=\mu_1$, and hence $W_1=W_0$. This complete our proof of Theorem 4.2.\qed

\begin{rem}
The argument in the last paragraph of the proof of Lemma 4.21 cannot work for $n=1$. But we do not know whether Lemma 4.21 is also true for $n=1$. On the other hand if we forget about the rectifiability condition, then there is a simple counter example: just notice that the map $f: S^1=[-\pi,\pi]\to [-1,1], f(\theta)=\sin(3\theta)$ satisfies that $\int_{-\frac\pi2}^{\frac\pi2}f(\theta+\phi)\cos(\phi)d\phi=0$ for all $\theta\in S^1$, hence the two conic 1-varifolds $V_1=\int_{S^1}V_{R_\theta} d\theta$ and $V_2=\int_{S^1}(1-f(\theta))V_{R_\theta} d\theta$ has the same weighted projection on all 1-linear subspaces.
\end{rem}

\section{Uniqueness of tangent varifolds for stationary 1-varifolds}

In this section we begin to prove unique tangent behaviour at an arbitrary point for stationary 1-varifolds. 

\begin{thm}[Uniqueness of tangent varifold] Let $V$ be a stationary 1-varifold in $\R^n$. Then for any $x_0\in [V]$, Var Tan$(V,x_0)$ contains only one element.
\end{thm}

\nd We will still use weighted projections, but only on 2 dimensional linear subspaces this time. As stated at the beginning of the proof of Theorem 4.2, the theorem holds for dimension $n=1,2$. So we are only going to prove it for ambient dimension at least 3.

So let $n\ge 2$. Let $V$ be any stationary 1- varifold in $\R^{n+1}$ such that ${\pi_P}_\s(V)$ is locally finite for all $P\in G(n+1,2)$, and $x_0\in[V]$. 

Let $W$ be any element in $\mbox{Var Tan}(V,x_0)$. It is a conical 1-varifold, hence there exists a finite Radon measure $\mu$ on $\partial B(0,1)$, such that
\be W=\int_{z\in S^n}V_{R_z}d\mu(z),\ee
Where for any point $z\in S^n$, $R_z$ is the half line issued from the origin and passing through the point $z$. Hence
\be |W|=\int_0^\infty \mu_rdr,\ee
where $\mu_r$ is a Radon measure supported on $rS^n$ with $\mu_r(A)=\mu(A/r)$.

Let $f$ be the shortest distance projection from $\R^{n+1}\bs\{0\}$ to $S^n$. Take any $P\in G(n+1,n)$, and let $v$ be a unit normal vector for $P$. For $r\in \R$, let $P_r=P+rv$, and for all $A\subset P$, denote by $A_r=A+rv$. Denote by $H$ the half space $\{x\in \R^{n+1}:<x,v>>0\}=\cup_{r>0}P_r$. For any $r>0$, $f:P_r\to f(P_r)=S^n\cap H$ is a continuous bijection. Let $\nu_r=\frac{|x|}{r}f^{-1}_\s(\mu)$, that is, for any $A\subset P_r$, $\nu_r(A)=\int_{x\in A}\frac{|x|}{r}d\mu(f^{-1}(x))$. 

Note that the measures $\nu_r$ are measures on planes parallel to $P$.   
Then for all $A\subset P$,
\be {\pi_P}_\s(\nu_r)(rA)={\pi_P}_\s(\nu_1)(A) ,\ee
and the restriction of $|W|$ on $H$ satisfies that 
\be |W|\lfloor_H=\int_{r>0}\nu_r.\ee

To see (5.5), take any Borel set $B\subset S^n$ and $0\le a<b<\infty$. Let $A_{B,a,b}=\{x\in H: f(x)\in B\mbox{ and }a< |x|< b\}$. It is then enough to verify that the two measures $|W|\lfloor_H$ and $\int_{r>0}\nu_r$ coincide on all such sets. We have
\be \begin{split}(\int_{r>0}\nu_r)(A_{B,a,b})&=\int_0^\infty \nu_r(A_{B,a,b}\cap P_r)dr
=\int_0^\infty dr[\int_{A_{B,a,b}\cap P_r}\frac{|x|}{r}df^{-1}_\s(\mu)]\\
&=\int_0^\infty dr[\int_{\{x:f(x)\in B, <x,v>=r,a<|x|<b\}}\frac{|x|}{r}df^{-1}_\s(\mu)]\\
&=\int_0^\infty dr[\int_{f\{x:f(x)\in B, <x,v>=r,a<|x|<b\}}\frac{|x|}{r}d\mu]\\
&=\int_0^\infty dr[\int_{\{\theta\in B: \frac{r}{<\theta,v>}\in (a,b)\}}\frac{1}{<\theta,v>}d\mu(\theta)]\\
&=\int_B \frac{1}{<\theta,v>}d\mu(\theta)[\int_0^\infty dr 1_{\{\frac{r}{<\theta,v>}\in (a,b)\}}]\\
&=\int_B \frac{1}{<\theta,v>}[(a-b)<\theta,v>]d\mu(\theta)=(a-b)\mu(B)=|W|(A_{B,a,b})
\end{split}\ee
by (5.3).
The third inequality is because $f$ is bijection on $P_r$ and $S^n\cap H$. 

Thus we have (5.5).

Next, we define $\gamma_1=\frac{1}{|x|^2}\nu_1$ on $P_1$. Then 
\be\begin{split}
\gamma_1(P_1)&\le \int_{P_1}|x|d\gamma_1=\int_{P_1}\frac{1}{|x|}d\nu_1=\int_{P_1}\frac{1}{|x|}|x|d\mu(f^{-1}(x))\\
&=\int_{P_1}d\mu(f^{-1}(x))=f_\s(\mu)(P_1)=\mu(S\cap H)<\infty.
\end{split}
\ee

Hence $\gamma_1$ is a finite measure on $P_1$, and $|x|$ is integrable with respect to $\gamma_1$.

Since $\R^{n+1}=\R v\oplus P$, points in $\R^{n+1}$ can be written in the form $(r,x_P)\in \R\times P$, where $x_P=\pi_P(x)\in P$ and $r=<x,v>$.  
Let $\gamma$ be the Radon measure on $P$  defined as $\gamma(A)=\gamma_1(\{(x_P,1):x_P\in A\})$. It is essentially a parallel version of $\gamma_1$, and $\int_P (1+|x|^2)^\frac 12d\gamma(x)<\infty$ by (5.7).

We want to calculate the Fourier transform of $\gamma$ on $P$. So fix any vector $\xi\in P\bs\{0\}$. Then the vectors $e_1=v,e_2=\frac{\xi}{|\xi|}$ are two orthonormal vectors. We complete it to an orthonormal basis $\{e_1,e_2,\cdots, e_{n+1}\}$ of $\R^{n+1}$. Then in this coordinate system, $P=\{x=(x_1,x_2,\cdots, x_{n+1}):\in \R^{n+1}, x_1=0\}$. We fix this coordinate system.

%

Set $Q=\{x\in \R^{n+1}: x_3=x_4=\cdots=x_{n+1}=0\}$ the 2-plane generated by $\xi$ and $v$. By our assumption, 
${\pi_Q}_\s(V)$ is locally finite, hence by Proposition 2.6, ${\pi_Q}_{**}(V)$ is a stationary varifold. Note that $W\in$ Var Tan$(V,x_0)$, hence Lemma 2.23 tells that ${\pi_Q}_{**}(W)\in $Var Tan$({\pi_Q}_{**}(V),\pi_Q(x_0))$, and hence is conic and locally finite.
%

Let $B\subset G_1(Q)$, then 
\be 
\begin{split}
{\pi_Q}_{**}(W)(B)&=\int_{z\in S^n}{\pi_Q}_{**}(V_{R_z})(B)d\mu(z)\\
&=\int_{z\in S^n}d\mu(z)\{\int_{\{(x,S):(\pi_Q(x), \pi_Q(S))\in B\}}|\pi_Q(z)|^2dV_{R_z}(x,S)\}\\
&=\int_{z\in S^n}d\mu(z)z_Q^2\{\int_{R_z}1_B(\pi_Q(x),\pi_Q(Z))dx\},
\end{split}
\ee
where $z_Q$ denotes $|\pi_Q(z)|$, and $Z\in G(n+1,1)$ denotes the linear subspace generated by $z$.

In particular, if $A\subset Q$ is Borel, then
\be 
\begin{split}
|{\pi_Q}_{**}(W)|(A)&={\pi_Q}_{**}(V)(A\times G_1(Q))\\
&=\int_{z\in S^n}d\mu(z)|z_Q|^2\{\int_{R_z}1_A(\pi_Q(x))dx\}\\
&=\int_{z\in S^n}d\mu(z)|z_Q|^2\H^1(\pi_Q^{-1}(A) \cap R_z)
\end{split}
\ee

Note that $\int_{z\in S^n}d\mu(z)\H^1\lfloor R_z=|W|$, hence
\be |{\pi_Q}_{**}(W)|(A)=\int_{\pi_Q^{-1}(A)}\frac{|z_Q|^2}{|z|^2}d|W|(z).\ee

Now for any $s<t$, set $A_{s,t}=\{z\in \R^{n+1}: \frac12\le z_1 \le \frac 32\mbox{ and }s\le \frac{z_2}{z_1}\le t\}$, and $B_{s,t}=A_{s,t}\cap Q $. Then $\pi_Q^{-1}B_{s,t}=A_{s,t}$, and hence
\be \begin{split}|{\pi_Q}_{**}(W)|(B_{s,t})&=\int_{A_{s,t}}\frac{|z_Q|^2}{|z|^2}d|W|(z)=\int_{\{z:\frac12\le z_1 \le \frac 32, s\le \frac{z_2}{z_1}\le t \}}\frac{|z_Q|^2}{|z|^2}d|W|(z)\\
&=\int_\frac12^\frac32 dz_1\int_{P_{z_1}}1_{A_{s,t}}\frac{|z_Q|^2}{|z|^2}d\nu_{z_1}(z_1,z_2,z_3,\cdots,z_{n+1})\\
&=\int_\frac12^\frac32 dz_1\int_{P_{z_1}\cap A_{s,t}}\frac{|z_Q|^2}{|z|^2}d\nu_{z_1}(z_1,z_2,z_3,\cdots,z_{n+1})
\end{split}
\ee

Note that $\pi_P(P_{z_1}\cap A_{s,t})=z_1\pi_P(P_1\cap A_{s,t})$, hence by (5.4), for a fixed $z_1$,
\be 
\begin{split}
&\int_{P_{z_1}\cap A_{s,t}}\frac{|z_Q|^2}{|z|^2}d\nu_{z_1}(z_1,z_2,z_3,\cdots,z_{n+1})\\
&=\int_{P_{z_1}\cap A_{s,t}}\frac{|z_Q/z_1|^2}{|z/z_1|^2}d\nu_{z_1}(z_1(1,\frac{z_2}{z_1},\frac{z_3}{z_1},\cdots,\frac{z_{n+1}}{z_1})\\
&=\int_{P_1\cap A_{s,t}}\frac{|y_Q|^2}{|y|^2}d\nu_1(1,y_2,y_3,\cdots, y_{n+1})\\
&=\int_{P_1\cap A_{s,t}}\frac{1+y_2^2}{|y|^2}d\nu_1(1,y_2,y_3,\cdots, y_{n+1})
\end{split}
\ee
and hence
\be 
\begin{split}
|{\pi_Q}_{**}(W)|(B_{s,t})&=\int_\frac12^\frac32 dz_1\int_{P_1\cap A_{s,t}}\frac{1+y_2^2}{|y|^2}d\nu_1(1,y_2,y_3,\cdots, y_{n+1})\\
&=\int_{P_1\cap A_{s,t}}\frac{1+y_2^2}{|y|^2}d\nu_1(1,y_2,y_3,\cdots, y_{n+1})\\
&=\int_{P_1\cap A_{s,t}}\frac{1+<y,\frac{\xi}{|\xi|}>^2}{|y|^2}d\nu_1(y)\\
&=\int_{P_1\cap A_{s,t}}1+<y,\frac{\xi}{|\xi|}>^2 d\gamma_1(y)\\
&=\int_{\{y\in P:s<<y,\frac{\xi}{|\xi|}><t\}}1+<y,\frac{\xi}{|\xi|}>^2d\gamma(y).
\end{split}\ee

Denote by $L_\xi\subset P$ the 1-subspace generated by $\xi$. Let $\gamma_\xi={\pi_{L_\xi}}_\s(\gamma)$ be the marginal of $\gamma$ on $L_\xi$. Note that for each $\lambda\in \R$, the function $1+<y,\frac{\xi}{|\xi|}>^2$ is constant and equals $(1+\lambda^2)$ on ${\pi_{L_\xi}}^{-1}(\lambda\frac{\xi}{|\xi|})$. Hence
\be \begin{split}
|{\pi_Q}_{**}(W)|(B_{s,t})&=\int_{\{y\in P:s<<y,\frac{\xi}{|\xi|}><t\}}1+<y,\frac{\xi}{|\xi|}>^2d\gamma(y)\\
&=\int_{s\frac{\xi}{|\xi|}}^{t\frac{\xi}{|\xi|}}(1+\lambda^2) d\gamma_\xi(\lambda\frac{\xi}{|\xi|}).
\end{split}\ee 

Recall that  $W$ is any element of Var Tan$(V,x_0)$. Since $V$ is a stationary varifold such that  ${\pi_Q}_{\s}(V)$ is locally finite for all $Q\in G(n+1,2)$, by Proposition 2.6, ${\pi_Q}_{**}(V)$ is a stationary 1-varifold in the 2-dimensional plane $Q$. But we already know that stationary 1-varifolds in $\R^2$ have unique tangent behaviour at every point, hence Var Tan$({\pi_Q}_{**}(V),\pi_Q(x_0))$ contains only one element. This yields
\be {\pi_Q}_{**}(W')={\pi_Q}_{**}(W),\forall Q\in G(n+1,2),\forall W'\in \mbox{Tan Var}(V,x_0).\ee

That is, all elements in Tan Var$(V,x_0)$ admit the same weighted projection on all 2-dimensional planes. We would like to show that this property guarantees that they are the same.

So let $W'$ be another element in Tan Var$(V,x_0)$. With respect to the same plane $P$ and its unit normal vector $v$, define the corresponding $\mu'$ on $S^n$,$\mu'_r$ on $rS^n$, $\nu'_r$ on $P_r$, $\gamma'_1$ on $P_1$, and $\gamma'$ on $P$ the same way as $\mu,\mu_r,\nu_r,\gamma_1$ and $\gamma$ for $W$.  For each $\xi\in P$, set $\gamma,_\xi={\pi_{L_\xi}}_\s(\gamma)$. Then similarly we have, for any $s<t$, 
\be |{\pi_Q}_{**}(W')|(B_{s,t})=\int_{s\frac{\xi}{|\xi|}}^{t\frac{\xi}{|\xi|}}(1+\lambda^2) d\gamma'_\xi(\lambda\frac{\xi}{|\xi|}).\ee

By (5.15), we know that ${\pi_Q}_{**}(W)={\pi_Q}_{**}(W')$. Hence (5.13) and (5.14) yields
\be  \int_{s\frac{\xi}{|\xi|}}^{t\frac{\xi}{|\xi|}}(1+\lambda^2) d\gamma_\xi(\lambda\frac{\xi}{|\xi|})=\int_{s\frac{\xi}{|\xi|}}^{t\frac{\xi}{|\xi|}}(1+\lambda^2) d\gamma'_\xi(\lambda\frac{\xi}{|\xi|})\ee
for any $\xi\in P$ and any $s<t$. In other words, for $\xi\in P$, and any segment $I\subset L_\xi$, we have
\be \int_I (1+|y|^2)d\gamma_\xi(y)=\int_I (1+|y|^2)d\gamma'_\xi(y).\ee

Note that $\gamma_\xi$ and $\gamma'_\xi$ are Radon measures on the line $L_\xi$, hence (5.18) yields that 
\be \gamma_\xi=\gamma'_\xi\ee
for all $\xi\in P\bs \{0\}$.

Thus, for the two finite measures $\gamma$ and $\gamma'$ on $P$, their Fourier transform satisfy that for all $\xi\in P\bs \{0\}$,
\be 
\begin{split}
\hat\gamma(\xi)&=\int_P e^{-i<x,\xi>}d\gamma(x)=\int_{L_\xi}e^{-i\lambda|\xi|}d\gamma_\xi(\lambda)\\
&=\int_{L_\xi}e^{-i\lambda|\xi|}d\gamma'_\xi(\lambda)=\int_P e^{-i<x,\xi>}d\gamma'(x)=\hat{\gamma'}(\xi),
\end{split}
\ee
which yields that $\gamma=\gamma'$ on $P$. A direct definition chasing gives that $|W|\lfloor_H=|W'|\lfloor_H$. Note that this is true for any half space $H$, and $W$ and $W'$ are conical, hence $|W|=|W'|$ on $\R^{n+1}$, and thus $W=W'$ on $\R^{n+1}$. This gives the uniqueness of tangent cones for 1-dimensional stationary varifold on $\R^{n+1}$.\qed

The following Theorem is just a local version of Theorem 5.1. It is easy to believe, since tangent behavior is only a local property. 

\begin{thm}Let $V$ be a 1-varifold in $\R^n$, and $V$ is stationary in an open subset $U\subset \R^n$, that is, for all $g\in C_C^\infty(U,\R^n)$, $\d V(g)=0$. Then for all $x_0\in [V]\cap U$, Var Tan$(V,x_0)$ contains only on element.
\end{thm}

\nd So let $V$ and $x_0$ be as in the statement. We will do the same cut and paste procedure as in Proposition 3.2. Take $r_0>0$ such that $B(x_0,r_0)\subset U$. By Lemma 3.3, take $r\in (0,r_0)$ such that $|\d V_r|$ is finite and supported on $\partial B(x_0,r)$, and (3.4) holds. Fix this $r$. For all $x\in \partial B(x_0,r)$, set $E_x=\{x+t\omega_r(x),t\ge 0\}$, which is the half line issued from the point $x$ with direction $\omega_r(x)$. By (3.5), $E_x\cap B(y,r)=\emptyset$. Define $V'$ and $W$ as in the beginning of the proof of Proposition 3.2. Then $W$ is a stationary varifold in $\R^n$, and Var Tan$(V,x_0)=$Var Tan$(W,x_0)$, which contains only one element by Theorem 5.1.\qed 

\section{Generalization to Riemannian manifolds}

The above theorem of uniqueness of tangent behaviour for 1-stationary varifolds in Euclidean spaces can be generalised without much difficulties to varifolds on general Riemannian manifold: we will prove Theorem 0.1 in this section. By Nash embedding theorem, we can restrict ourselves to sub Riemannian manifolds of Euclidean spaces. So we give directly the definitions on these manifolds. For definitions on a abstract Riemannian manifold, see for example \cite{AlAl}.

\begin{defn}[Varifolds] 

$1^\circ$ Let $M$ be a Riemannian submanifold of dimension $n$ in $\R^N$. A $k$-dimensional varifold $V$ on $M$ is just a Radon measure on $G_k(M)$, where $G_k(M)\subset G_k(\R^N)$ is the bundle $\cup_{x\in M}G(T_xM,k)$. Denote by $|V|$ the image of $V$ under the bundle projection $\pi:G_k(M)\to M$. Then $|V|$ is a Radon measure on $M$.

$2^\circ$ Let $V$ be a 1-dimensional varifold on $M$. Denote by $\mathfrak{X}(M)$ the vector space of smooth mappings $g:M\to TM$ with compact support such that $g(x)\in T_x(M)$ whenever $x\in M$.

Then the first variation $\d V$ of $V$ is a linear map from $\mathfrak{X}(M)\to \R$ defined as follows: 
\be \d V(g)=\int_M (D_sg(x)\cdot s) dV(x,S), \ee
where $s$ is a unit vector in $S$, and $D$ is the covariant differentiation with respect to the Levi-Civita connection on $M$.

We say that $V$ is stationary if $\d V=0$.

$3^\circ$ Suppose $V$ is a varifold on $M$, and $x\in M$. Let $j: T_xM\to T_x\R^N=\R^N$ and $i:M\to \R^N$ be the two inclusion map. Then a varifold $C$ on $T_xM$ is said to be a tangent varifold of $V$ at $x$ if $j_\s(C)$ is a tangent varifold of $i_\s V$. Denote by Var Tan$(V,x)$ the set of all tangent varifolds of $V$ at $x$.
\end{defn}

\begin{rem} $1^\circ$ In the rest of this section, we will also use the letter $V$ to denote the varifold $i_\s(V)$ on $\R^N$, and we do not distinguish $C$ and $i_\s(C)$ either, to save notations.

$2^\circ$ Again, as the case in $\R^n$, the non emptiness of Var Tan$(V,x)$ at any point $x\in [V]$ is simply guaranteed by functional analysis.

$3^\circ$ On an abstract Riemannian manifold, the tangent varifold is defined via the exponential map. That is, $C$ is a tangent varifold of $V$ at a point $x$ if there exists a sequence $\{\lambda_i\}$ that tends to zero such that $C=\lim_{i\to\infty}(\lambda_i^{-1}\log_x)_\s V$. Note that the varifolds
 ${\eta_{x,\lambda}}_\s V$ and $(\lambda^{-1}\log_x)_\s V$ are not the same. But the two define the same limits, i.e. same sets of tangent varifolds, since the differential of $\log_x$ at the point $x$ is the identity map.\end{rem}

Now let us get back to the tangent behaviour for stationary 1-varifolds.

\begin{thm}Let $M\subset\R^N$ be a $n$-dimensional Riemannian submanifold of $\R^N$, and let $V$ be a stationary 1-varifold on $M$. Then for any $x_0\in [V]$, Var Tan$(V,x_0)$ contains only one element. 
\end{thm}

\nd Let $M$, $V$ and $x_0$ as stated in the theorem. Denote by $\d^MV$ the first variation of $V$ as a varifold in $M$, and $\d^N V$ the first variation of $V$ as a varifold in $\R^N$. Denote by $D^M$ the Levi-Civita connection on $M$, and $D^N$ that on $\R^N$ (that is, the flat one).

Let us calculate $\d^N V$. So let $g\in C_C^\infty(\R^N,\R^N)$. Then at each point $x\in M$, we can decompose $g(x)=g^T(x)+g^\perp(x)$, where $g^T(x)\in T_xM$, and $g^\perp(x)\in T_xM^\perp$. Then both $g^T$ and $g^\perp$ are in $C_C^\infty(\R^N,\R^N)$.

For the part $g^T$, we have
\be \begin{split}
\d^NV(g^T)&=\int_{\R^N} (D^N_sg^T(x)\cdot s) dV(x,S)=\int_{\R^N} ([D^M_sg^T(x)+II(s,g^T(x))]\cdot s) dV(x,S)\\
&=\d^M V(g^T)+\int_{\R^N}(II(s,g^T(x))\cdot s) dV(x,S), \end{split}\ee
where $II_x(\cdot,\cdot):T_xM\times T_xM\to T_x^\perp M$ is the vector valued second fundamental form on $M$, and $s$ is a unit vector that generates $S$. As a result, $II(s,g(x))$ is a normal vector to $M$ on $x$, hence $II(s,g^T(x))\cdot s=0$. Note that $\d^M V=0$ by hypothesis, therefore
\be \d^NV(g^T)=0, \forall g\in C_C^\infty(\R^N,\R^N);\ee

On the other hand, for $g^\perp$,
\be \begin{split}\d^NV(g^\perp)&=\int_{\R^N} (D^N_sg(x)\cdot s) dV(x,S)\\
&=\int_{\R^N} ((D_s^Ng^\perp(x))^T\cdot s))dV(x,S)\\
&=\int_{\R^N} (II(s,s)\cdot g^\perp)dV(x,S).\\
\end{split}\ee
Hence
\be \d^NV(g^\perp)=\int_{\R^N} (II(s,s)\cdot g^\perp) dV(x,S),\ \forall g\in C_C^\infty(\R^N,\R^N).\ee

Altogether we have for any $g\in C_C^\infty(\R^N,\R^N)$, 
\be \d^NV(g)=\int_{\R^N} (II(s,s)\cdot g^\perp) dV(x,S),\ee
and 
\be |\d^NV(g)|\le |V|(\mbox{spt }g)||g||_\infty ||II\lfloor_{\mbox{spt }g}||_\infty.\ee

As a result, $|\d^N V|$ is absolutely continuous to $|V|$. In particular, it is locally finite and hence is a Radon measure. By Riesz representation theorem (cf. \cite{Fe} 2.5.12), there exists a measurable $\omega:\R^N\to S^{N-1}$, such that
\be \d^NV(g)=\int_{\R^N}<g(x)\cdot \omega(x)>d|\d^NV|(x).\ee

It is not hard to see that $|\d^NV|$ is supported on $M$, and by (6.9), for $|\d^NV|(x)$-a.e. $x\in M$, $\omega(x)$ is a unit normal vector to $M$. For these $x$, set $E_x=\{x+t\omega(x),t\ge 0\}$ the half line issued from the point $x$ with direction $\omega(x)$.

Now fix $x_0\in [V]$. Then there exists $r_0>0$ such that for $|\d^NV|(x)$-a.e. $x\in B(x_0,r_0)\cap M$, $E_x\cap (M\cap B(x_0,r_0)=\emptyset$, and for all $x\in B(x_0,r_0)$, $\pi_{T_{x_0}M}(\omega(x))\le \frac12$, where $\pi_{T_{x_0}M}:\R^N\to T_{x_0}M$ is the orthogonal projection. Such an $r_0$ exists since $M$ is a Riemannian submanifold of $\R^N$.

Set 
\be V_0=\int_{M\cap B(x_0,r_0)}V_{E_x}d|\d^NV|(x),\ee

 and $W=V_0+V$. Then $W$ is a varifold in $\R^N$, and it is stationary in $B(x_0,\frac12 r_0)$. By Proposition 1.26, the 1-density $\theta_1(|W|,x_0)$ exists. By definition, $|W|=|V|+|V_0|\ge |V|$. Hence the 1-density $\theta_1(|W|,x_0)\ge \theta_1(|V|,x_0)>0$, since $x\in [V]$. On the other hand, by (6.10), since $|\d^N V|$ locally controlled by $|V|$, hence so is $|V_0|$ and is $|W|$. Therefore $\theta_1(|W|,x_0)<\infty$ since $x\in [v]$ implies that $\theta_1(|V|,x_0)<\infty$. as a result, $x_0\in [W]$. By Theorem 5.21, Var Tan$(W,x_0)$ contains only one element, we call it $W_0$.
%

We want to prove that Var Tan$(V,x_0)$ contains only one element, which is the part of $W_0$ that is tangential to $T_{x_0}M$. So let $Z\in $Var Tan$(V,x_0)$. Then $Z$ is a varifold in $T_{x_0}M$, and there exists $r_k\to 0$ such that $Z=\lim{\eta_{x_0,r_k}}_\s V$. Recall that by (6.10) $|V_0|$ is locally controlled by a multiple of $|V|$, hence it has finite upper density at the point $x_0$. As a result, modulo taking a subsequence, we can suppose that the sequence ${\eta_{x_0,r_k}}_\s V_0$ converges to a varifold $Z_0$ as well. Then $Z+Z_0$ is a tangent varifold of $W$, hence $Z+Z_0=W_0$.

On the other hand, denote by $H=\{S\in G(N,1): |\pi_{T_{x_0}M}s|>\frac 12\}$, where for any element $S\in G(N,1)$, $s$ stands for a unit normal vector in $S$. Then notice that for any $R>0$, and for $k$ large such that $R/r_k<r_0$, by definition of $r_0$, we have that ${\eta_{x_0,r_k}}_\s V_0(B(x_0,R)\times H)=0$. As a result, the limit $Z_0$ satisfies that
\be Z_0 (\R^N\times H)=0.\ee

In other words, $W_0\lfloor{\R^N\times H}=Z\lfloor{\R^N\times H}$. But $Z$ is a varifold on $T_{x_0}M$, which is contained in $H$, hence it is supported on $\R^N\times H$. As a result, $Z=Z\lfloor{\R^N\times H}=W_0\lfloor{\R^N\times H}$. This yields that Var Tan$(V,x_0)=\{W_0\lfloor{\R^N\times H}\}$, which contains only one element. \qed

\renewcommand\refname{References}
\bibliographystyle{plain}
\bibliography{reference}

\end{document}